\documentclass[journal,twocolumn]{IEEEtran}

\usepackage{amsfonts}
\usepackage{amsmath}

\usepackage{bm}
\usepackage{pgf}
\usepackage{psfrag}
\usepackage{enumerate}
\usepackage{amssymb}
\usepackage{bbm}
\usepackage{color}
\usepackage{algorithm}
\usepackage{algpseudocode}
\usepackage{subfigure}
\usepackage{url} 
\def\ie{{\it i.e.}}
\def\eg{{\it e.g.}}
\algnewcommand\algorithmicinput{\textbf{Input:}}
\algnewcommand\INPUT{\item[\algorithmicinput]}
\algnewcommand\Offline{\item[\textbf{Offline-Phase:}]}
\algnewcommand\Online{\item[\textbf{Online-Phase:}]}
\def\nn{\nonumber}
\def\expec{\mathbbm{E}}
\def\real{\mathbbm{R}}

\DeclareMathOperator{\dom}{dom}

\DeclareMathOperator*{\argmin}{argmin}
\def\defeq{\triangleq}

\newcommand{\st}{\textup{subject to}}
\DeclareMathOperator{\minimize}{\text{\textnormal{minimize}}}

\newtheorem{remark}{Remark}
\newtheorem{definition}{Definition}
\newtheorem{assumption}{Assumption}
\newtheorem{theorem}{Theorem}
\newtheorem{lemma}[theorem]{Lemma}

\newcommand{\vv}[1] {{#1}}
\def\b{\vv{b}}
\def\u{\vv{u}}

\def\f{\vv{f}}

\def\p{\vv{p}}
\def\d{\vv{d}}
\def\h{\vv{h}}

\def\r{\vv{r}}

\def\b{s}
\def\u{u}

\def\f{f}
\def\r{r}

\def\bmax{S^{\max}}
\def\bmin{S^{\min}}
\def\bdot{S^{(\cdot)}}
\def\umax{U^{\max}}
\def\umin{U^{\min}}
\def\udot{U^{(\cdot)}}

\def\p{p}

\def\d{\delta}
\def\la{\lambda}
\def\muC{\mu^\mathrm{C}}
\def\muD{\mu^\mathrm{D}}

\newcommand{\pos}[1]{\left(#1\right)^+}
\renewcommand{\neg}[1]{\left(#1\right)^-} 
\def\minimize{\mbox{minimize  }}
\def\st{\mbox{subject to  }}

\def\h{h}
\def\hC{\h^\mathrm{C}}
\def\hD{\h^\mathrm{D}}

\def\g{g}

\newcommand{\opttag}[1]{\mbox{\bf#1  }}

\def\fn{\vv{\f}}

\def\bsn{\vv{\bs}}

\def\ustatn{\vv{\ustat}}
\def\rstatn{\vv{\rstat}}
\def\thetastatn{\vv{\thetastat}}
\def\fstatn{\vv{\fstat}}
\def\uhatn{\vv{\uhat}}
\def\fhatn{\vv{\fhat}}
\def\thetahatn{\vv{\thetahat}}
\def\rhatn{\vv{\rhat}}

\def\fmax{F^{\max}}

\def\pmin{P^{\min}}
\def\pmax{P^{\max}}
\def\J{J}
\newcommand{\Jk}[1]{\J_\mathrm{P#1}}

\newcommand{\Pk}[1]{\text{\bf P#1}}

\def\ustat{\u^\mathrm{stat}}
\def\fstat{\f^\mathrm{stat}}
\def\rstat{\r^\mathrm{stat}}
\def\thetastat{\theta^\mathrm{stat}}

\def\bs{\tilde \b}
\def\ks{\Gamma}
\def\ksmin{\ks^{\min}}
\def\ksmax{\ks^{\max}}
\def\W{W}
\def\Wmax{\W^{\max}}
\def\uhat{\u^\mathrm{ol}}
\def\rhat{\r^\mathrm{ol}}
\def\thetahat{\theta^\mathrm{ol}}
\def\fhat{\f^\mathrm{ol}}
\def\Dl{\underline{D}}
\def\Du{\overline{D}}
\def\M{M}
\def\Mone{\M^{\u}}
\def\Mtwo{\M^{\b}}

\def\N{\eta}
\def\None{\N^{\u}}
\def\Ntwo{\N^{\b}}
\def\L{L}
\def\Ls{\Delta}

\def\dmin{\d^{\min}}
\def\dmax{\d^{\max}}
\def\X{X}
\def\Xumin{\X^{\min,\u}}
\def\Xumax{\X^{\max,\u}}
\def\Xbmin{\X^{\min,\b}}
\def\Xbmax{\X^{\max,\b}}
\def\Xudot{\X^{(\cdot),\u}}

\def\Xbdot{\X^{(\cdot),\b}}

\def\Tday{\mathcal{T}^\mathrm{Day}}

\def\Fmax{F^\mathrm{max}}
\def\rpos{\r^+}
\def\rneg{\r^-}
\def\gpos{\g^+}
\def\gneg{\g^-}
\def\Sm{\mathbb{S}}
\newcommand{\Tn}[1]{\mathcal{T}^{\mathrm{#1},k}}
\def\Fcal{\mathcal{F}}

\usepackage{scalefnt}

\begin{document}

\title{Distributed Online Modified Greedy Algorithm\\ for Networked Storage Operation under Uncertainty}
\author{Junjie~Qin, 
        Yinlam~Chow, 
        Jiyan~Yang, 
        and~Ram~Rajagopal
\thanks{This research was supported in part by the Satre Family fellowship, and in part by the Tomkat Center for Sustainable Energy.}
\thanks{J. Qin, Y. Chow and J. Yang are with the Institute for Computational and Mathematical Engineering, Stanford University, Stanford, CA, 94305 USA, (e-mail: \{jqin,ychow,jiyan\}@stanford.edu).}%
\thanks{R. Rajagopal is with the Department of Civil and Environmental Engineering, Stanford University, Stanford, CA 94305 USA (e-mail: ramr@stanford.edu).}%
}

\maketitle

\begin{abstract}
The integration of intermittent and stochastic renewable energy resources requires increased flexibility in the operation of the electric grid. Storage, broadly speaking, provides the flexibility of shifting energy over time; network, on the other hand,  provides the flexibility of shifting energy over geographical locations. The optimal control of storage networks in stochastic environments is an important open problem. The key challenge is that, even in small networks, the corresponding constrained stochastic control problems on continuous spaces suffer from curses of dimensionality, and are intractable in general settings. For large networks, no efficient algorithm is known to give optimal or provably near-optimal performance for this problem. This paper provides an efficient algorithm to solve this problem with performance guarantees. We study the operation of storage networks, \ie,  a storage system interconnected via a power network. An online algorithm, termed Online Modified Greedy algorithm, is developed for the corresponding constrained stochastic control problem. A sub-optimality bound for the algorithm is derived, and a semidefinite program is constructed to minimize the bound. In many cases, the bound approaches zero so that the algorithm is near-optimal. A task-based distributed implementation of the online algorithm relying only on local information and neighbor communication is then developed based on the alternating direction method of multipliers. Numerical examples verify the established theoretical performance bounds, and demonstrate the scalability of the algorithm.
\end{abstract}

\addtolength{\topskip}{-.04in}
\addtolength{\abovedisplayskip}{-.04in}
\addtolength{\belowdisplayskip}{-.04in}

\section{Introduction}
Deep penetration of renewable energy generation is essential to ensure a sustainable future. Renewable energy resources, such as wind and solar, are intrinsically variable. Uncertainties associated with these intermittent and volatile resources pose a significant challenge to their integration into the existing grid infrastructure \cite{NRELWest2010}. More flexibility, especially in shifting energy supply and/or demand across time and network, is desired to cope with the increased uncertainties.

Energy storage provides the functionality of shifting energy across time. A vast array of technologies, such as  batteries, flywheels, pumped-hydro, and compressed air energy storages, are available for such a purpose \cite{Denholm2010, lindley2010naturenews}.
Furthermore, flexible or controllable demand provides another ubiquitous source of storage. Deferrable loads -- including many thermal loads, loads of internet data-centers and loads corresponding to charging electric vehicles (EVs) over certain time interval \cite{thermalStor1993} -- can be interpreted and controlled as \emph{storage of demand} \cite{ObRACC2013}. 
Other controllable loads which can possibly be shifted to an earlier or later time, such as thermostatically controlled loads (TCLs), may be modeled and controlled as a storage with negative lower bound and positive upper bound on the storage level \cite{HaoSanandajiPoollaVincent2013}.
These forms of storage enable inter-temporal shifting of excess energy supply and/or demand, and significantly reduce the reserve requirement and thus system costs. 

On the other hand, shifting energy across a network, \ie, moving excess energy supply to meet unfulfilled demand  among different geographical locations with transmission or distribution lines, can achieve similar effects in reducing the reserve requirement for the system.
Thus in practice, it is natural to consider these two effects together. Yet, it remains mathematically challenging to formulate a sound and tractable problem that accounts for these effects in electric grid operations. Specifically, due to the power flow and network constraints, control variables in connected buses are coupled. Due to the storage constraints, control variables in different time periods are coupled as well.
On top of that, uncertainties associated with stochastic generation and demand dramatically complicate the problem, because of the large number of recourse stages and the need to account for all probable realizations.

Two categories of approaches have been proposed in the literature. The first category is based on exploiting structures of specific problem instances, usually using dynamic programming. These structural results are valuable in providing insights about the system, and often lead to analytical solution of these problem instances.
However, such approaches rely heavily on specific assumptions of the type of storage, the form of the cost function, and the distribution of uncertain parameters. Generalizing results to other specifications and more complex settings is usually difficult, and consequently this approach is mostly used to analyze single storage systems.
For instance, analytical solutions to optimal storage arbitrage with stochastic price have been derived in  \cite{QRsimpleStorPes2012} without storage ramping constraints, and in \cite{MITrampStor} with ramping constraints. Problems of using energy storage to minimize energy imbalance are studied in various contexts; see \cite{SuEGTPS, RLDSACC} for reducing reserve energy requirements in power system dispatch, \cite{BitarRACC_colocated, Powell} for operating storage co-located with a wind farm, \cite{IBMload, DataCenter} for operating storage co-located with end-user demands, and \cite{StorDRLongbo} for storage with demand response.
The other category is to use heuristic algorithms, such as Model Predictive Control (MPC) \cite{XieEtAlWindStorMPC} and look-ahead policies \cite{NRELStorValue2013},  to identify sub-optimal storage control rules. Usually based on deterministic (convex) optimization, these approaches can be easily applied to general networks. The major drawback is that these approaches usually do not have any performance guarantee. Consequently,  it lacks theoretical justification for their implementation in real systems.
Examples of this category can be found in \cite{XieEtAlWindStorMPC} and references therein.



This work aims  at designing distributed online deterministic optimizations that solve the stochastic control problem with provable guarantees. It contributes to the existing literature in the following ways. First, we formulate the problem of  storage network operation as a stochastic control problem with general cost functions, which encapsulates a variety of problems with different types of storage as well as different uses of storage.
Second, we devise an online algorithm for the problem based on the theory of Lyapunov optimization, 
and prove guarantees for its performance in terms of a bound of its sub-optimality. This converts the ``intractable'' stochastic control program to a sequence of tractable deterministic optimization programs.
The bound is useful not only in assessing the performance of our algorithm, but also in evaluating the performance of other sub-optimal algorithms when the optimal costs are hard to obtain. It can also be used to estimate the maximum cost reduction that can be achieved by {\it any} storage operation, thus provides understanding for the limit of a certain storage system.
To the best of our knowledge, this is the first algorithm with provable guarantees for the storage operation problem with general electric networks. 
Finally, we derive task-based distributed implementation of the online algorithm using the alternating direction method of multipliers (ADMM). 

The paper generalizes our prior work \cite{OMGarxiv} by modeling a networked storage system, and extending the online control algorithm to the network setting. Preliminary results related to the network setting have been presented in \cite{QCYR:acm}. The online optimization in this paper is different from that in \cite{QCYR:acm}, and the sub-optimality bound here is significant superior to the bound for the algorithm proposed in  \cite{QCYR:acm}. The aspect of distributed implementation is also new in this paper. 


The rest of the paper is organized as follows. Section 2 formulates the problem of operating a storage network under uncertainty. Section 3 gives the online algorithm and states the performance guarantee. Section 4 discusses the distributed implementation of the online program. Numerical examples are then given in Section 5. Section 6 concludes the paper.

\section{Problem Formulation}
\def\Pcal{\mathcal{P}}
\subsection{Centralized Problem}

We model the power grid as a directed graph $G(V,E)$, with $V = [n] \defeq \{1, \dots, n\}$, $E = [m] \defeq \{1,\dots, m\}$, where $n$ is the number of nodes and $m$ is the number of edges. The node-edge incidence matrix $A\in \real^{n\times m}$ defined as
\[
A_{i,e} = \begin{cases}
 1 & \mbox{if } e \to i, \\
 -1 & \mbox{if } e \leftarrow i, \\
 0 & \mbox{otherwise},
 \end{cases}
\]
where $e \to i$ denotes that $i$ is the head of $e$, and $e \leftarrow i$ denotes that $i$ is the tail of $e$.\footnote{Notation: For a directed graph $G(V,E)$, define $V(e) \defeq \{i\in V: i \sim e\}$, and $E(i) \defeq \{e\in E: e\sim i\}$, where $i \sim e$ (and $e\sim i$) means that edge $e$ and node $i$ are incident. We assume that all these sets are equipped with the natural order. 
For any vector $v\in \real^d$ and $\Pcal \subseteq [d]$, $v_{\Pcal} \in \real^{|\Pcal|}$ is the sub-vector containing entries of $v$ indexed by set $\Pcal$. Similarly, for any matrix $M \in \real^{d_1 \times d_2}$, and $\Pcal_1 \subseteq [d_1]$ and $\Pcal_2 \subseteq [d_2]$, $M_{\Pcal_1, \Pcal_2} \in \real^{|\Pcal_1|\times |\Pcal_2|}$ is the sub-matrix containing rows and columns of $M$ indexed by sets $\Pcal_1$ and $\Pcal_2$. For any variable $x\in \real^m$ that is defined for each edge, and if edge $e \in E$ is incident to nodes $i$ and $j$, we use the notations $x_e$ and $x_{ij}$ interchangeably to refer to the $e$th element of $x$. For any $x\in \real$, $\pos{x} \defeq \max(x,0)$ and $\neg{x} \defeq \pos{-x}$. An extended real value function $g(x)$ with domain $\dom g = \mathcal{C} \subseteq \real^d$ is such that $g(x) = \infty$ if $x\not\in \mathcal{C}$. } 
Here each node models a bus and each edge models a line. To simplify the exposition, we assume each bus $i$ is connected to all of the following types of devices:\footnote{By setting the problem data properly, we can model buses which are only connected to a subset of these devices. For example, a generator bus with no renewables and no storage can be modeled by setting $\d_i(t) =0$ and $\bmin_i=\bmax_i =0$.  }
\begin{itemize}

\item Uncontrollable net supply. A renewable generator and a load are connected to the bus, with the net power supply, \ie, the generation minus the demand, at time period $t$ denoted by $\d_i(t)$. As both demand and generation can be stochastic, $\d_i(t)$ is in general stochastic.
\item Energy storage. A storage with storage capacity $\bmax_i$, minimum storage level $\bmin_i$, storage charging limit $\umax_i$, and storage discharging limit $-\umin_i$ is connected to the bus. The storage level (or state of charge) is denoted by $\b_i(t)$ and the storage control is denoted by $\u_i(t)$ with $\u_i(t) > 0$ representing charging and $\u_i(t)< 0$ representing discharging. For each time period $t$, we have constraints $\bmin_i \le \b_i(t) \le \bmax_i$ and $\umin_i \le \u_i(t) \le \umax_i$. The storage dynamics is $\b_i (t+1) = \la_i \b_i(t) + \u_i(t)$, where $\la_i \in (0,1]$ is the storage efficiency which models the energy loss over time without storage operation. We denote the set of parameters for the storage at bus $i$ by $\Sm_i \defeq \{\la_i, \bmin_i, \bmax_i, \umin_i, \umax_i\}$. Here the set of parameters for each storage satisfies feasibility  and controllability assumptions (see Assumption~\ref{assume:feas} in Appendix~\ref{app:lyap_pf} and \cite{OMGarxiv} for more discussions). 
\item Conventional generator. Its generation at time period $t$ is denoted by $\rpos_i(t)$ ($\ge 0$) and its convex cost function is denoted by $\gpos_i(\cdot)$. 
It is possible in certain scenarios to have more supply than demand (\eg, when there is too much wind generation). In such cases, let $\rneg_i(t)$ ($\ge 0$) be the generation curtailment at time period $t$ and $\gneg_i(\cdot)$ be the cost associated with the curtailment.  Without loss of optimality, we can summarize $\rpos_i(t)$ and $\rneg_i(t)$ by a single variable $\r_i(t)$ such that $\rpos_i(t) = \pos{\r_i(t)}$ and $\rneg_i(t) = \neg{\r_i(t)}$. Then the total cost at bus $i$ and in time period $t$ is 
\[
\g_i(\r_i(t)) = \gpos_i(\rpos_i(t)) + \gneg_i(\rneg_i(t)).
\]
Optionally, the cost can depend on a stochastic price parameter $\p_i(t) \in [\pmin_i, \pmax_i]$, so that we write the cost as $\g_i(\r_i(t); \p_i(t))$. 
   
\end{itemize}
We use the classic DC approximation for AC power flow. For time period $t$, let the voltage phase angle on bus $i$ be $\theta_i(t)$. Then the real power flow from bus $i$ to bus $j$ can be written as
\[
\f_{ij}(t) = B_{ij} (\theta_i(t) - \theta_j(t)),
\]
where $B\in \real^{n\times n}$ is the imaginary part of the admittance matrix (Y-bus matrix) under DC assumptions, and $\f_{ij}(t)$ satisfies line flow constraints $-\fmax_{ij} \le \f_{ij}(t) \le \fmax_{ij}$, where $\fmax_{ij} = \fmax_{ji} > 0$ is the real power flow capacity of the line connecting bus $i$ and bus $j$. Note that our focus here is to identify the optimal operation of storage systems under uncertainty. A more detailed modeling of the AC power flow and incorporating recent convexification techniques \cite{slowOPF1, slowOPF2} into our algorithm are left for future work.

We can now formulate the problem as a stochastic control problem as follows
\begin{subequations}\label{p1}
\begin{align}
\minimize &  (1/T)\expec \sum_{t=1}^T \sum_{i=1}^n \g_{i}(\r_i(t); \p_i(t)) \\
\st & \delta_i(t) +\r_i(t)  = \u_i(t) + \sum_{j=1}^n  \f_{ij}(t) , \label{n1}\\
& \b_i(t+1) = \la_i\b_i(t) + \u_i(t),  \label{n2}\\
& \umin_i \le \u_i(t) \le \umax_i, \label{n3}\\
& \bmin_i \le \b_i(t) \le \bmax_i, \label{n4}\\
&B_{ij} (\theta_i(t) -\theta_j(t)) = \f_{ij}(t) ,\label{e1}\\
& -\Fmax_{ij} \le \f_{ij}(t) \le \Fmax_{ij} \label{e2},
\end{align}
\end{subequations}
where $T$ is the total number of time periods under consideration, the expectation is taken over $\p_i(t)$ and $\d_i(t)$, constraints \eqref{n1}, \eqref{n2}, \eqref{n3} and \eqref{n4} hold for all $i$ and $t$,  constraints \eqref{e1} and \eqref{e2} hold for all $i$, $j$ and $t$, and $\b_i(1) \in [\bmin_i, \bmax_i]$ is given for each $i$. Here the goal is to find an optimal control policy for each time period $t$ which maps the information available up to the time period to the optimal decisions $( \u^\star(t),\r^\star(t), \theta^\star(t),\f^\star(t))$. 

Albeit the bulk of this paper focuses on the formulation~\eqref{p1}, we note that it can be extended in various directions.
\begin{remark}[Generalized Storage Model]\label{rk:gs}
The storage model described above consider primarily energy storage. But following the development in \cite{OMGarxiv}, it is easy to incorporate other type of generalized storage such as deferrable loads as storage of demand, and collections of thermostatically controlled load. 
In addition, the  energy loss during charging/discharging can be modeled with conversion functions. For example, a storage with charging coefficient $\muC\in (0,1]$ and discharging coefficient $\muD \in (0,1]$ can be modeled using charging conversion function  $\hC(\u) = (1/\muC)\u$ and discharging conversion function $\hD(\u) = \muD \u$, respectively. See \cite{OMGarxiv}  for more details. 
\end{remark}
\begin{remark}[Nonconvex Objective]
The assumption that $\g_i$ is convex for each $i \in [n]$ is not strictly necessary.
See \cite{QCYR:report} for generalization to general subdifferentiable functions.
\end{remark}
\begin{remark}[Other Costs and Constraints]
Many other costs including operational cost of storage due to charging and discharging, and other constraints including bounds on the generation and phase angles can be added without altering our results and the proofs. In fact, the cost can be a function of the form $g_i(u_i(t), r_i(t),\theta_i(t), \delta_i(t), p_i(t))$.
\end{remark}

Our prior work \cite{OMGarxiv} can be viewed as the single bus special case of the problem formulated here. Thus the examples for different use cases of the storage (\eg, balancing and arbitrage) discussed in \cite{OMGarxiv} can also be encapsulated into our current framework together with a network. The incorporation of the network element allows our methodology to be applied to broader problems such as microgrid management and storage-based real-time regulation for the bulk power grid.




\subsection{Cluster based Distributed Control}\label{sec:cluster}
Solving problem~\eqref{p1} in a centralized fashion may not be feasible due to concerns regarding privacy, communication, and computation. First of all, specifying the centralized problem~\eqref{p1} requires collection of information about the cost functions and parameters of the devices connected to each of the buses, and the probability distributions of all local stochastic parameters. This process involves agents who own the generators, storages, as well as power consumers who may not be willing to report such data. Even if the data reporting is granted, 
gathering all these data from nodes of a large power network, and subsequently disseminating the optimal control signal obtained from the centralized solution in real time presents a challenge on the communication system required. The large amount of data that have to be sent to and from the centralized control center may lead to traffic congestions and delays in the data delivery. Finally, granting an adequate communication infrastructure in place, solving the stochastic control problem formulated in~\eqref{p1} over a large network is not tractable due to a lack of practical algorithms, \ie, existing algorithms either do not have any performance guarantee or do not scale gracefully with the number of buses of the system. 

A cluster-based control architecture for the future grid is envisioned in \cite{grip:cdc2011}. Here we present a first step in achieving such an architecture. In particular, we consider solving the centralized problem~\eqref{p1} with resource clusters. Suppose that the network is partitioned into $L$ clusters. Each cluster $C_\ell$ consists of a subset of nodes $V_\ell \subset V$ and a subset of lines $E_\ell \subset E$, \ie, $C_\ell\defeq (V_\ell, E_\ell)$, and is controlled by a cluster controller (CC).  The CC for each cluster $C_\ell$  
\begin{itemize}
\item possesses local static information including $\g_i$ and $\Sm_i$ for all $i \in V_\ell$, and $B_e$ and $\fmax_e$ for all $e\in E_\ell$,
\item senses local disturbances $\d_i(t)$ and $\p_i(t)$ for all $i \in V_\ell$ and all $t$,
\item controls local variables $\u_i(t)$, $\r_i(t)$ and $\theta_i(t)$ for all $i \in V_\ell$, and $\f_e(t)$ for all $e\in E_\ell$ and all $t$,
\item and communicates with its neighbors $C_{N_\ell}$ where $C_{N_\ell}$ is the collection of $C_{w}$'s for which there exists $e\in E_\ell$, $i \in V_w \mbox{ such that } e\sim i$, or there exists $e\in E_w$, $i \in V_\ell \mbox{ such that } e\sim i$.
\end{itemize}

Here we provide a bird-eye view of our approach for tackling the challenging distributed stochastic control problem which we just formulated. Section~\ref{sec:omgnet} provides an online algorithm that converts the centralized stochastic control program to a sequence of online deterministic optimization. Section~\ref{sec:admm} then presents the decentralization of these online deterministic optimization using the alternating direction method of multipliers (ADMM). 

\section{Online Modified Greedy Algorithm for Networked Storage Control}\label{sec:omgnet}
\subsection{Algorithm}
We propose a very simple algorithm to solve the centralized problem~\eqref{p1} with performance guarantees.  The algorithm, termed the network online modified greedy (OMG) algorithm, is composed of an offline and online phase. Next we describe the input data to the algorithm and each phase. \vspace{.2cm}

\noindent{\bf Input Data.}
Similar to the single storage online modified greedy (OMG) algorithm \cite{OMGarxiv}, for each bus $i \in [n]$, in addition to the storage parameters $\Sm_i$ and the cost functional form $\g_i$,  the algorithm requires two input parameters that are a lower bound, denoted by $\Dl \g_i$, and an upper bound, denoted by $\Du \g_i$, for the subdifferential of the objective function $\g_i$ with respect to $\u_i(t)$.\footnote{Mathematical expressions for these parameters are relegated to Appendix~\ref{app:def4omg}. \label{fn3}} 

\begin{remark}[Distribution-Free Method] As in the single storage case \cite{OMGarxiv}, The OMG algorithm is a distribution-free method in the sense that almost no information regarding the joint probability distribution of the stochastic parameters $\d_i(t)$ and $\p_i(t)$ are required. The only exception is when calculating $\Du\g_i$ and $\Dl \g_i$, the support of $\p_i(t)$ and $\delta_i(t)$ may be needed. Comparing to the entire distribution functions, it is much easier to estimate the supports of the stochastic parameters  from historical data.
\end{remark}\vspace{.2cm}

\noindent{\bf Offline Phase.} Before running the algorithm, each bus $i\in [n]$ needs to calculate two algorithmic parameters, namely a shift parameter $\ks_i$ and a weight parameter $\W_i$. Any pair $(\ks_i, \W_i)$ satisfies the following conditions can be used:
\begin{align}
\ksmin_i\le & \ks_i \le \ksmax_i, \label{eq:ksbounds}\\
0 < & \W_i  \le \Wmax_i, \label{eq:Wbounds}
\end{align}
where $\ksmin_i$, $\ksmax_i$ and $\Wmax_i$ are functions of the storage parameters $\Sm_i$ and subdifferential bounds  $\Dl\g_i$ and $\Du \g_i$.

It will be clear later that the sub-optimality bound depends on the choice of $(\ks_i, \W_i)$. As in \cite{OMGarxiv}, we provide two approaches for selecting these parameters
\begin{itemize}
\item The \emph{maximum weight} approach (\texttt{maxW}): Setting $\W_i = \Wmax_i$ reduces the interval in \eqref{eq:ksbounds} to a singleton ($\ksmin_i = \ksmax_i$) and hence determines a unique $\ks_i$. 

\item The \emph{minimum sub-optimality bound} approach (\texttt{minS}): It turns out that the sub-optimality bound of OMG, as a function of $(\ks_i, \W_i)$'s for all $i\in [n]$, can be minimized using a semidefinite program reformulation. This approach uses the set of $(\ks_i, \W_i)$'s minimizing the sub-optimality bound.
\end{itemize}\vspace{.2cm}

\noindent{\bf Online Phase.} At the beginning of each time period $t$, the OMG algorithm solves a deterministic optimization as follows
\begin{subequations}\label{OMGcen}
\begin{align}
 \minimize  &  \sum_{i=1}^n (\la_i/\W_i) (\b_i+\ks_i) \u_i + \g_{i}(\r_i; \p_i)\\
\st  
        & \umin_i \le u_i \le \umax_i,\\
        & \delta_i +\r_i  = \u_i + \sum_{j=1}^n  \f_{ij},\\
        &B_{ji} (\theta_j -\theta_i) = \f_{ji} ,\\
        & -\Fmax_{ji} \le \f_{ji} \le \Fmax_{ji}.
\end{align}
\end{subequations}
where the optimization variables are $\u$, $\r$, $\theta$ and $\f$, and we have dropped the dependence on $t$ to simplify the notation. 
This treatment is justified by the fact that \eqref{OMGcen} does not involve the charging and discharging constraints induced by the storage capacity and storage dynamics, \ie, we have removed constraints \eqref{n2} and \eqref{n4}, which can be alternatively summarized as
\begin{equation}\label{n2+4}
 \bmin_i \le \la_i \b_i + \u_i\le \bmax_i.
\end{equation}
It will be show later in Appendix~\ref{app:lyap_pf} that \eqref{n2+4} holds automatically given that the algorithmic parameters of OMG satisfy conditions in \eqref{eq:ksbounds} and \eqref{eq:Wbounds}. 
 
The optimization is similar to the greedy heuristics which minimize the stagewise cost, \ie, $\sum_{i=1}^n  \g_{i}(\r_i; \p_i)$, subject to constraints of \eqref{OMGcen} together with constraint \eqref{n2+4} in each step. 
Instead of directly optimizing the cost at the current time period, for each storage, the OMG algorithm optimizes a weighted combination of the stage-wise cost and a linear term of $\u_i$ depending on the shifted storage level $\b_i +\ks_i$. Here the weight parameter $W_i$ decides the importance of the original cost in this weighted combination, while the shift parameter $\ks_i$ defines the shifted state given the original state $\b_i$.  
Roughly speaking, the shifted state $\b_i+\ks_i$ belongs to an interval $[\bmin_i+\ks_i, \bmax_i+\ks_i]$ which usually contains $0$. For fixed $W_i>0$,
if the storage level is relatively high, the shifted state is greater than $0$, such that the state-dependent term (\ie, $(\la_i/W_i) (\b_i + \ks_i) \u_i$) encourages a negative $\u_i$ (discharge) to minimize the weighted sum. As a result, the storage level in the next time period will be brought down. On the other hand, if the storage level is relatively low, the shifted state is smaller than $0$, such that the state-dependent term encourages a positive $\u_i$ (charge) and consequently the next stage storage level is increased. These two effects together help to hedge against uncertainty by maintaining a storage level somewhere in the middle of the feasible interval.
More detailed discussion regarding the design of the modification term in the objective can be found  in \cite{OMGarxiv}.

\subsection{Performance Guarantees}\label{sec:pg}
We provide a stylized analysis for the performance of OMG.
\begin{assumption}\label{ass:1}
The following assumptions are in force for the analysis in this section.
\begin{enumerate}[{\bf A}1]
\item Infinite horizon: The horizon length $T$ approaches to infinity.
\item IID disturbance: The disturbance process $\{(\d(t), \p(t))\in \real^{2n}: t\ge 1\}$ is independent and identically distributed (i.i.d.) across $t$ and is supported on a compact set such that $\d_i(t) \in [\dmin_i, \dmax_i]$ and $\p_i(t) \in [\pmin_i, \pmax_i]$ for all $i\in [n]$ and all $t$. Note that any correlation structure is allowed for variables in the same time period.
\item Frequent acting: The storage parameters satisfy $\umax_i-\umin_i < \bmax_i -\bmin_i$ for all $i \in [n]$.
\end{enumerate}
\end{assumption}
Here {\bf A}1 and {\bf A}2 are technical assumptions introduced to simplify the exposition. Relaxing {\bf A}1 leads to no change in our results except an extra term of $O(1/T)$ in the sub-optimality bound. For large $T$, this term is negligible. \cite{OMGarxiv} discusses how to reduce {\bf A}2. Under these two assumptions, the storage operation problem can be cast as an infinite horizon average cost stochastic optimal control problem in the following form
\begin{subequations}\label{p1_inf}
\begin{align}
\minimize \quad &  \lim_{T \to \infty}  (1/T) \expec \Big[ \sum_{t=1}^T\sum_{i=1}^n \g_{i}(\r_i(t); \p_i(t)) \Big]\\
\st \quad &\eqref{n1}, \eqref{n2}, \eqref{n3}, \eqref{n4}, \eqref{e1}, \eqref{e2}.
\end{align}
\end{subequations}

Assumption {\bf A}3  states that the range of feasible storage control $\umax_i -\umin_i$ is smaller than the range of storage levels $\bmax_i -\bmin_i$, \ie, the ramping limits of the storage is relatively small compared to the storage capacity. 
For any storage system, this assumption is true as long as the length of each time period $\Delta t$ is made small enough; see \cite{OMGarxiv} for more details.

Define $J(u,r,\theta,f)$ as the  total cost of problem~\eqref{p1} induced by the sequence of control $\{(\u(t),\r(t),\theta(t),f(t)),\,t\geq 1\}$ and $J^\star=J(u^\star,r^\star,\theta^\star, f^\star)$ as the minimum cost of the average cost stochastic control problem with $\{(u^\star(t),r^\star(t),\theta^\star(t),f^\star(t)),\,t\geq 1\}$ being the corresponding optimal control sequence. 
The main results regarding the performance of the OMG algorithm is summarized as follows, whose proof is relegated to Appendix~\ref{app:lyap_pf}.
\begin{theorem} [Performance]\label{thm:perf_lyap}
The control sequence $(\uhat,\rhat,\thetahat, \fhat) \defeq\{(\uhat(t),\rhat(t),\thetahat(t), \fhat(t)), t\ge 1\}$ generated by the OMG algorithm is feasible with respect to all constraints of problem~\eqref{p1} and its sub-optimality  is bounded by $\sum_{i=1}^n\M_i(\ks_i)/\W_i$, that is
\begin{equation}\label{eq:1busiidperf_bdd}
J^\star \le  J(\uhat,\rhat,\thetahat, \fhat) \le J^\star + \sum_{i=1}^n\M_i(\ks_i)/\W_i,
\end{equation}
where
\begin{align*}
&\M_i(\ks_i)=\Mone_i(\ks_i)+\la_i(1-\la_i)\Mtwo_i(\ks_i),\\
&\Mone_i(\ks_i) =\! \frac{1}{2}\! \max\left(\! \left(\umin_i\!+\!(1\!-\!\la_i)\ks_i\right)^2\!\!,\left(\umax_i\!+\!(1\!-\!\la_i)\ks_i\right)^2\! \right)\!,\\
&\Mtwo_i(\ks_i)= \max\left( \left(\bmin_i+\ks_i\right)^2,\left(\bmax_i+\ks_i\right)^2 \right).
\end{align*}
\end{theorem}
The theorem above guarantees that the cost of the OMG algorithm is no greater than $J^\star + \sum_{i=1}^n\M_i(\ks_i)/\W_i$. 

In many cases, we are interested to minimize the sub-optimality bound. This can be cast as the following optimization
\begin{subequations}\label{P3:PO}
\begin{align*}
\opttag{PO:}\minimize & \quad \sum_{i=1}^n\M_i(\ks_i)/\W_i\\
\st & \quad  \ksmin_i\le \ks_i \le \ksmax_i,\,\, 0< \W_i \le \Wmax_i,
\end{align*}
\end{subequations}
where the constraints hold for all $i\in [n]$. Observing that the objective is separable across buses, we can solve this program separately on each bus via a semidefinite program (SDP) as in the single storage case \cite{OMGarxiv}. Here the SDP is reproduced for completeness.
\begin{lemma}[Semidefinite Reformulation of {\bf PO}]\label{SDP_P3_PO}
For each $i \in [n]$, let symmetric positive definite matrices $\Xumin_i$, $\Xumax_i$, $\Xbmin_i$ and $\Xbmax_i$ be defined as follows
\begin{equation*}
\!\!\Xudot_i\!\! = \!\!\begin{bmatrix}
\None_i&\!\!\udot_i+(1-\la_i)\ks_i\\
*& 2\W_i
\end{bmatrix}, \,\,
\Xbdot_i\!\! = \!\!\begin{bmatrix}
\Ntwo_i&\!\!\bdot_i+\ks_i\\
*& \W_i
\end{bmatrix},\!\!
\end{equation*}
where $(\cdot)$ can be either $\max$ or $\min$, and $\None$ and $\Ntwo$ are auxilliary variables. Then {\bf PO} can be solved via the following semidefinite program
\begin{subequations}\label{prob:sdp}
\begin{align}
\!\!\emph{\minimize} \quad & \None_i+\la_i(1-\la_i)\Ntwo_i \!\!\\
\!\!\emph{\st} \quad & \ksmin_i\le  \ks_i \le \ksmax_i,\,\, 0 <  \W_i  \le \Wmax_i, \!\!\\
& \Xumin_i,\Xumax_i, \Xbmin_i, \Xbmax_i \succeq 0, \!\!
\end{align}
\end{subequations}
where  $\ksmin_i$ and $\ksmax_i$ are linear functions of $\W_i$ as defined in \eqref{eq:ineq_1} and \eqref{eq:ineq_2}.
\end{lemma}

We close this section by summarizing some of the properties for the sub-optimality bound at each bus $i$ in the next remark; more detailed discussion and examples of the uses of the sub-optimality bound can be found at \cite{OMGarxiv}. 
\begin{remark}[Properties of $\M_i(\ks_i)/\W_i$]
The following properties are true for the per bus sub-optimality $\M_i(\ks_i)/\W_i$:
\begin{itemize}
\item For ideal storage ($\la_i = 1$), $\M_i(\ks_i)/\W_i$ is minimized with $\W_i = \Wmax_i$. 
\item Let the bound minimizing parameter choice be $(\ks_i^\star, \W_i^\star)$. Then $\M_i(\ks_i^\star)/\W_i^\star \to 0$ if (i) $\bmax_i - \bmin_i \to \infty$ while $\umax_i - \umin_i$ is fixed
or (ii) $\umax_i-\umin_i \to 0$ while $\bmax_i-\bmin_i$ is fixed (which may be the case when the storage is controlled frequently such that the length of each time period $\Delta t \to 0$).
That is, when the storage capacity is much larger than the range of feasible storage control action, the algorithm is optimal.
\end{itemize}

\end{remark}

\section{Distributed Online Control Via Alternating Direction Method of Multipliers}\label{sec:admm}

Results in previous section convert the stochastic control program~\eqref{p1} to a sequence of online deterministic optimization programs. In this section, we take a bottom-up approach in deriving a decentralized solution to~\eqref{p1}. In particular, we first reformulate the online program and then apply ADMM to obtain a \emph{fully distributed algorithm} that specifies computation and communication tasks for each bus and each line of the network. We then associate the corresponding tasks to the CC's to which these buses or lines belong. For a survey of ADMM, see \cite{BoydADMM}.
 
\subsection{Node-Edge Reformulation}\label{sec:admm:n+e}
In order to obtain a fully distributed algorithm that uses only local computation and neighborhood communication, it is necessary to ensure that all constraints of the optimization program only couple variables controlled by pairs of neighboring node and edge so that all communication can be implemented using simple pairwise messages. To this end, we reformulate the online program~\eqref{OMGcen} by creating local copies of certain variables. In particular, 
let $x_i\defeq (u_i, r_i, \theta_i, \widehat{f}_{i,E(i)})^\intercal$ be the local (primal) variables at node $i$, and $z_e \defeq (f_e, \widehat{\theta}_{e,V(e)})^\intercal$ be the local (primal) variables at edge $e$, where $\widehat{f}_{i,E(i)} \in \real^{|E(i)|}$ is node $i$'s local auxiliary copy of edge variable $f_{E(i)}$, and $\widehat{\theta}_{e, V(e)}\in \real^{2}$ is edge $e$'s local auxiliary copy of node variable $\theta_{V(e)}$. Here we use the notation $\widehat \f_{i,e}$ for $e\in E(i)$ to refer to $i$'s local copy of variable $\f_e$; similar notation $\widehat \theta_{e,i}$ is also used.  Then program~\eqref{OMGcen} can be written as 
\begin{subequations}\label{online:ref}
\begin{align}
\minimize \quad & \sum_{i=1}^n q_i(x_i) + \sum_{e=1}^m h_e(z_e)\\
\st \quad & \widehat{f}_{i, E(i)} = f_{E(i)}, \forall i \in [n],\label{online:ref:c1}\\
& \widehat{\theta}_{e, V(e)} = \theta_{V(e)}, \forall e\in [m],\label{online:ref:c2}
\end{align}
\end{subequations}
where extended real value functions $q_i$ and $h_e$ summarize the separable objective and constraints at node $i$ and edge $e$, respectively, and are defined as follows
\[
q_i(x_i)\defeq q_i(u_i,r_i,\theta_i, \widehat{f}_{i,E(i)}) \defeq (\la_i/\W_i) (\b_i+ \ks_i) \u_i + \g_{i}(\r_i; \p_i),
\]
with domain 
$\dom q_i = \{x_i: \umin_i\le \u_i \le \umax_i, \, 
\d_i +\r_i + A_{i,E(i)} \widehat{f}_{i, E(i)} = \u_i,\, \theta_i\in \real \}$, 
and $h_e(z_e) = 0$ with $ \dom h_e \defeq \{z_e: f_e = B_e A_{V(e), e}^T 
\widehat{\theta}_{e, V(e)}, \,\, -\fmax_e \le f_e \le \fmax_e\}$. Here constraints~\eqref{online:ref:c1} and~\eqref{online:ref:c2} ensures that at the solution, these local auxiliary variables must be equal to the corresponding true variables. The (scaled) dual variables\footnote{See Appendix~\ref{app:admm} for more details.} corresponding to constraints~\eqref{online:ref:c1} and~\eqref{online:ref:c2} are denoted by $\eta_i$ and $\xi_e$, respectively. We proceed to state the task-based distributed ADMM. The derivation of the algorithm is relegated 
to Appendix~\ref{app:admm}.

At each iterate, indexed by $k$, the following tasks are issued and completed in order:

\begin{itemize}
\item $\Tn{NP}_i$: Each node $i\in [n]$ performs \emph{node primal update}:
\begin{align*}
x^{k+1}_i &= \argmin_{x_i} q_i(x_i) + \frac{\rho}{2} \|\widehat f_{i, E(i)} - f_{E(i)}^k + \eta_i^k\|_2^2\\
&\quad\quad \quad\quad\quad\quad \quad + \sum_{e\in E(i)} \frac{\rho}{2}(\widehat \theta_{e,i}^k- \theta_i + \xi_{e,i}^k)^2,
\end{align*}
and then passes a message containing $\theta_i^{k+1}$ and $\widehat f_{i,e}^{k+1}$ to each neighboring edge $e\in E(i)$.
\item $\Tn{EP}_e$: Each edge $e\in [m]$ performs \emph{edge primal update}:
\begin{align*}
z^{k+1}_e &= \argmin_{z_e} h_e(z_e) +\frac{\rho}{2}\|\widehat \theta_{e, V(e)} - \theta_{V(e)}^{k+1} + \xi_e^k\|_2^2\\
&\quad\quad \quad\quad\quad\quad \quad + \sum_{i\in V(e)} \frac{\rho}{2} (\widehat f_{i,e}^{k+1}-f_e +\eta_{i,e}^{k})^2,
\end{align*}
and then passes a message containing $f_e^{k+1}$ and $\widehat \theta_{e,i}$ to each neighboring node $i\in V(e)$. 
\item $\Tn{ND}_i$: Each node $i\in [n]$ performs \emph{node dual update}:
\[
\eta^{k+1}_i  = \eta^k_i + \widehat f_{i, E(i)}^{k+1} - f_{E(i)}^{k+1},
\]
and passes a message containing $\eta^{k+1}_{i,e}$ to each neighboring edge $e\in E(i)$.
\item $\Tn{ED}_e$: Each edge $e\in [m]$ performs \emph{edge dual update}:
\[
\xi^{k+1}_e = \xi^k_e + \widehat \theta_{e, V(e)}^{k+1} - \theta_{V(e)}^{k+1},
\]
and passes a message containing $\xi^{k+1}_{e,i}$ to each neighboring node $i\in V(e)$. 
\end{itemize}
We summarize the convergence property of the iterates specified above, whose proof is relegated to Appendix~\ref{app:admm}.
\begin{lemma}\label{lemma:admmconv}
The iterates $(x^k, z^k)$ produced by tasks $\mathcal{T}^k = \left\{\Tn{NP}_{[n]}, \Tn{EP}_{[m]}, \Tn{ND}_{[n]}, \Tn{ED}_{[m]}\right\}$ are convergent. Let $x^\star \defeq \lim_{k\to \infty} x^k$ and $z^\star \defeq \lim_{k\to \infty} z^k$. Then $(x^\star, z^\star)$ is primal feasible and achieves the minimum cost of problem~\eqref{online:ref}. Furthermore, the rate of convergence is $O(1/k)$. 
\end{lemma}
\begin{remark}
Minimum amount of assumptions are required to obtain the convergence results given in Lemma~\ref{lemma:admmconv}. In particular, we do not assume the objective function is strongly convex which is a necessary assumption for standard distributed algorithms based on primal or dual decomposition. Furthermore, the rate of convergence for our algorithm is superior to primal or dual decomposition based algorithms, which usually have a rate of convergence $O(1/\sqrt{k})$. 
\end{remark}
\begin{remark}[Asynchronous Variant]
Based on the analysis in \cite{wei20131}, one can easily extend the algorithm described above to its asynchronous counterpart with similar convergence guarantees.  
\end{remark}

\subsection{Cluster-based Implementation}
In a cluster-based distributed control environment, each CC is responsible for a subset of resources in the grid. It is not necessary the case that there is a CC for each node and each edge. However, issuing tasks defined for each node and edge to the associated CC would implement our distributed algorithm in a cluster-based control environment.
 The iterates now have the following form: in order, each CC $\ell \in [L]$ (i) performs $\Tn{NP}_i$ for all $i \in V_\ell$, (ii)  performs $\Tn{EP}_e$ for all $e\in E_\ell$, (iii) performs $\Tn{ND}_i$ for all $i \in V_\ell$, and (iv)  performs $\Tn{ED}_e$ for all $e\in E_\ell$. Note that if the source and destination of a message belong to different CCs, instead of direct communications between the node-edge pair, the message is sent from the CC containing the source to the CC containing the destination\footnote{Recall the setup in Section~\ref{sec:cluster}: each CC $\ell$ can communicate with its neighbors $C_{N_\ell}$ where $C_{N_\ell}$ is the collection of $C_{w}$'s for which there exists $e\in E_\ell$, $i \in V_w \mbox{ such that } e\sim i$, or there exists $e\in E_w$, $i \in V_\ell \mbox{ such that } e\sim i$. As all messaging tasks only involve incident node-edge pairs, the communication between these CCs are possible. }; if a single CC controls both the source and destination of a message, the corresponding messaging step may be skipped.

%

\section{Numerical Tests}
In this section, we show three sets of numerical tests with different focuses. The first example (Subsection~\ref{sec:nt:star}) uses synthetic data that honor the i.i.d. assumption in Section~\ref{sec:pg} to demonstrate the use of the online algorithm and to show how the sub-optimality bound scales with storage parameters. The second example (Subsection~\ref{sec:nt:14}) applies the algorithm on IEEE 14 bus network together with real demand and wind data. The i.i.d. assumption no longer holds in this setup. We also demonstrate the convergence of ADMM in this setting. The last example (Subsection~\ref{sec:nt:scale}) is constructed in particular to show how the distributed algorithm scales with the number of buses of the system. All examples are implemented and tested using Matlab 2014a on a workstation with AMD Magny Cours 24-Core 2.1 GHz CPU  and 96GB RAM.

\subsection{Star Network}\label{sec:nt:star}
\begin{figure}[htbp]
\centering
\subfigure[Time homogeneous costs]{ \label{fig:result_Nstor_perc_a}
\centering
\psfrag{SMR}[cc][Bl]{Total $\bmax$ in network}
\psfrag{PercentageSaving}[cc][Bt]{Percentage cost savings (\%)}
\psfrag{NoStorage}{\scriptsize No storage}
\psfrag{Lyapunov}{\scriptsize OMG}
\psfrag{Greedy}{\scriptsize Greedy}
\psfrag{LowerBoundLegend}{\scriptsize Lower bound}
\psfrag{UpperBoundLegend}{\scriptsize Upper bound}
\includegraphics[width=.22\textwidth]{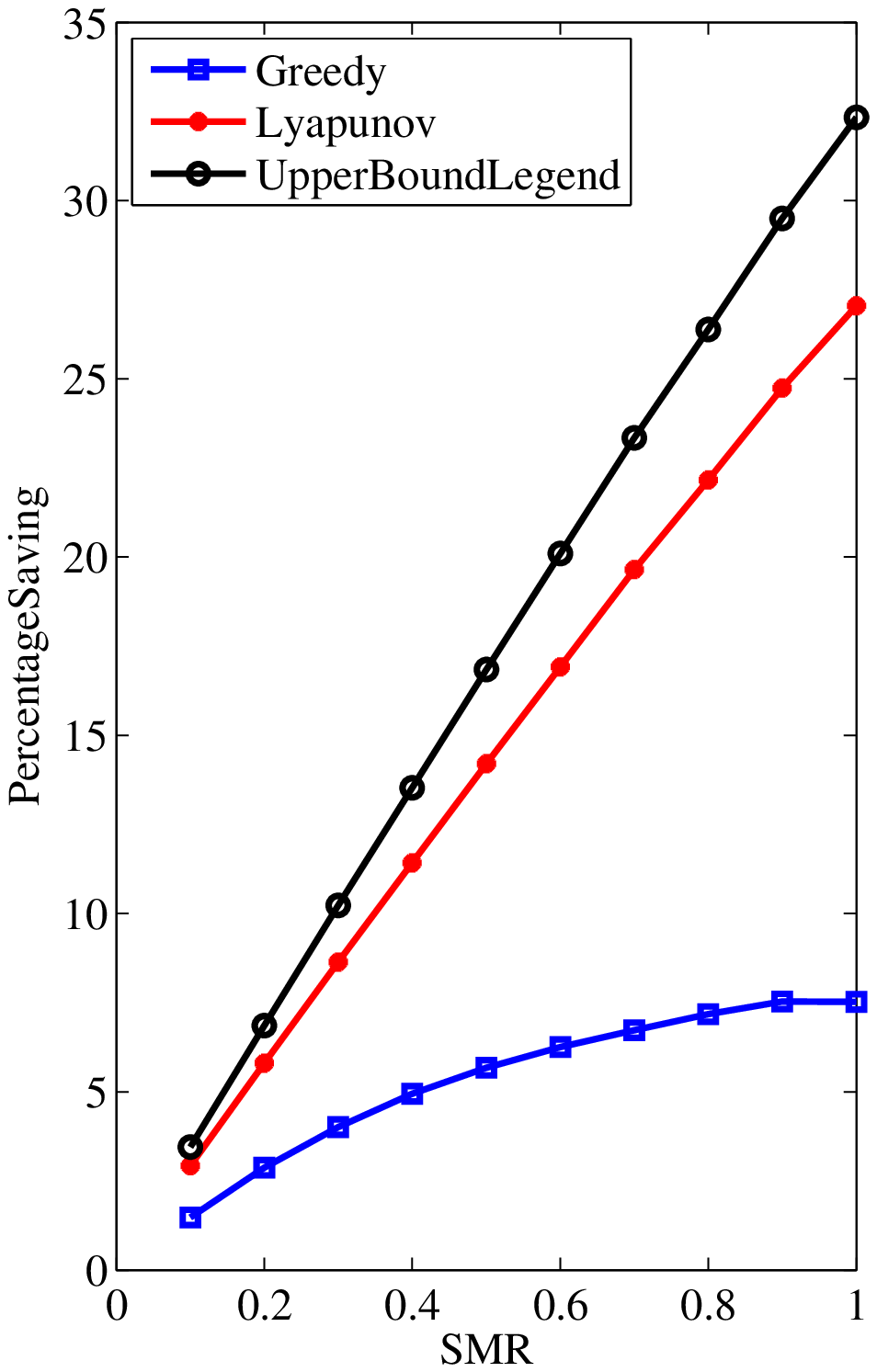}}
\subfigure[Time inhomogeneous costs]{ \label{fig:result_Nstor_perc_b}
\centering
\psfrag{SMR}[cc][Bl]{Total $\bmax$ in network}
\psfrag{PercentageSaving}[cc][Bt]{Percentage cost savings (\%)}
\psfrag{NoStorage}{\scriptsize No storage}
\psfrag{Lyapunov}{\scriptsize OMG}
\psfrag{Greedy}{\scriptsize Greedy}
\psfrag{LowerBoundLegend}{\scriptsize Lower bound}
\psfrag{UpperBoundLegend}{\scriptsize Upper bound}
\includegraphics[width=.22\textwidth]{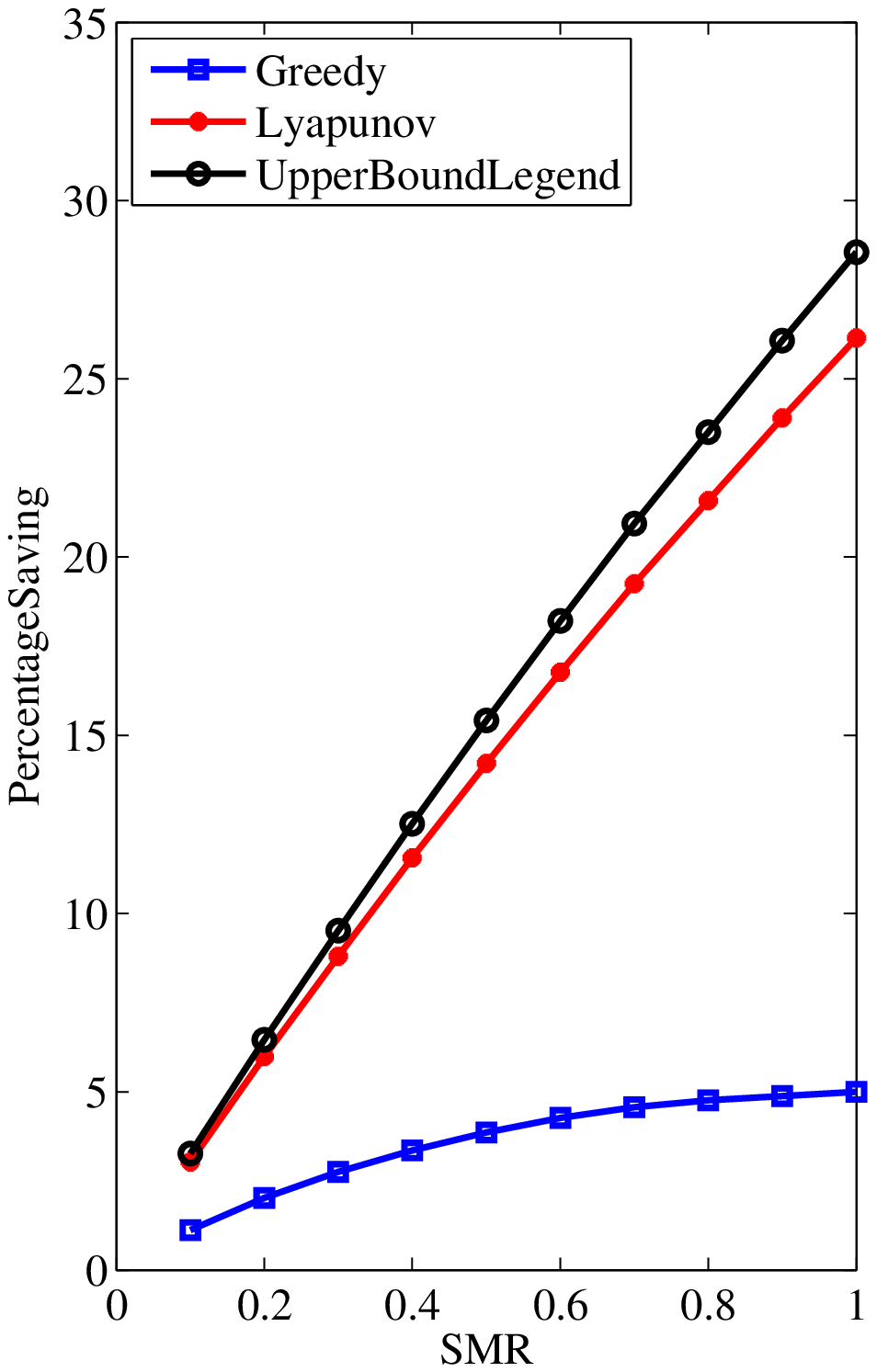}}
\caption{Percentage cost savings of a storage network operated for balancing. }\label{fig:5}
\end{figure}

Consider a star network, \ie, a tree with a root node and $(n-1)$ leaf nodes. With a homogeneous setting, all nodes are connected to identical power system components, and thus we only provide specification for a single bus $i$.
The storage network is operated for the purpose of balancing the demand and supply residual due to forecast errors in the wind power generation. The motivation of this setting in a single storage scenario is discussed in detail in \cite{SuEGTPS}. Let $\d_i(t)$ models the wind forecast error process for each bus $i$. We simulate the $\d_i(t)$ processes by generating Laplace distributed random variables with zero mean and standard deviation $0.149$ p.u. as in \cite{SuEGTPS}, which are estimated empirically using the NREL dataset. 
Two cases with different cost functions are considered. In the first case, time homogeneous costs of the form
\begin{equation} \label{eq:cost:net:homo}
\g_i(\r_i(t);\p_i(t)) = \p^\mathrm{H}_i(t)\neg{\r_i(t)},
\end{equation}
are considered, where $\p^\mathrm{H}_i(t) \equiv 1$;
in the second case, the cost function is modified to has a higher penalty rate during the day
\begin{equation*}
\g_i (\r_i(t)) = \p_i(t) \neg{\r_i(t)}=
\begin{cases}
3 \neg{\r_i(t)}  &t\in \Tday,  \\
\neg{\r_i(t)} , & \mbox{otherwise},
\end{cases}
\end{equation*}
with $\Tday$ is the set of time points during the day (7am to 7pm in our tests), $\p_i(t) =3$ if $t\in \Tday$, and $\p_i(t) =1$ otherwise.
We consider non-idealized storages which are operated frequently such that $\la_i = 0.999$ with conversion coefficients being $\muC_i = \muD_i = 0.995$ (\emph{cf}., Remark~\ref{rk:gs}). We fix $-\umin_i = \umax_i = (1/10) \bmax_i$. We have $n= 5$ and $\fmax_e = \sigma_\d$ for each line $e\in [m]$. The time horizon for the simulation is chosen to be $T=1000$. Figure~\ref{fig:5} shows the percentage cost savings compared to the no storage scenario.
Albeit the greedy heuristics have been proved to be the optimal solution for single storage systems in the time homogeneous cost setting in \cite{SuEGTPS}, OMG outperforms the greedy heuristics in the case with a network. The improvement over the greedy cost is more significant for the time inhomogeneous case. For both cases, the costs of OMG are close to the upper bounds derived using the sub-optimality bounds of the algorithm.


\subsection{IEEE 14 Bus Case}\label{sec:nt:14}
\begin{figure}[htbp]
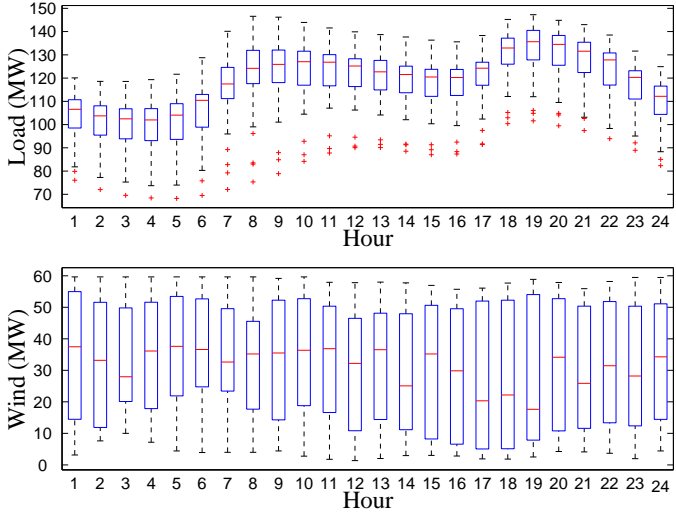

\begin{center}
\scalebox{0.31}{\hspace{-0.6in}\scalefont{1.5} 
\input{load_plot.pgf}}
\scalebox{0.31}{\hspace{-0.6in}\scalefont{1.5} 
\input{wind_plot.pgf}}
\end{center}
\caption{Bar plots for scaled hourly total load (upper panel) and wind data (lower panel) used for the simulation.
}\label{fig:data}
\end{figure}

The network data from IEEE 14 bus test system \cite{IEEEtest} are used for this example, with modifications described as follows. Three generators are connected to the network, \ie,  a coal power plant with capacity $500$MW and (constant) marginal generation cost $50$\$/MWh connected to bus 1, a nuclear power plant with capacity $450$MW and marginal generation cost $25$\$/MWh connected to bus 2, and a natural gas power plant with capacity $400$MW and marginal generation cost $100$\$/MWh is connected to bus 8.\footnote{The labeling of the buses are consistent with \cite{IEEEtest}} 
A wind power plant is connected to bus 3. Hourly data of wind power generation for January 2004 (Figure~\ref{fig:data}) are obtained from the NREL dataset \cite{NREL2010}, and are scaled to model a 30\% penetration scenario.
The hourly load data are obtained from PJM interconnection for the same period (Figure~\ref{fig:data}), and are scaled down and then factored out according to the portion of different load buses. 
Three storages are connected to buses 6, 7 and 10. Their capacities are $\bmax_6 = 300$MWh, $\bmax_7 =240$MWh, $\bmax_{10} = 300$MWh, and charging/discharging power rating are $\umax_6 = \umax_7= \umax_{10}=10$MW with $\umin_i = -\umax_i$ for all $i$. For simplicity (and the fact that conversion to cluster based implementation is easy), we emulate a complete distributed setting, where each node or each edge solves its own tasks in the distributed ADMM algorithm. 

The performance of OMG together with the greedy heuristic are simulated over $T=744$ time periods (\ie, hourly for January 2004). We also compute the cost if there is no storage in the system, and the offline clairvoyant optimal cost which solves the storage operation problem assuming the full knowledge of the future load and wind ahead of time. 
For this example, the hourly average no storage cost is \$$51710$. The costs of the greedy heuristics, OMG, and offline optimal are $96.1$\%, $95.7$\% and $90.3$\% of the no storage cost, respectively. Here the cost achieved by the offline optimal solution is a loose lower bound as it requires information that is not available to the decision maker. The stochastic lower bound, estimated by our algorithm under i.i.d. assumption is $94.6$\% of the no storage cost. As the disturbances are not i.i.d., we expected the actual optimal cost is between these two lower bounds. 

\begin{figure}[htbp]
\centering
\subfigure[Objective function convergence]{ \label{fig:ADMM_subgrad_conv}
\centering
\psfrag{Legend1ADMM10}{\scriptsize ADMM: $\rho=10$}
\psfrag{Legend2ADMM100}{\scriptsize ADMM: $\rho=100$}
\psfrag{Legend3ADMM500}{\scriptsize ADMM: $\rho=500$}
\psfrag{Legend4SubGrad1}{\scriptsize SubGD: $\zeta^k=10^{-3}/k$}
\psfrag{Legend5SubGrad05}{\scriptsize SubGD: $\zeta^k=10^{-3}/k^{0.5}$}
\psfrag{Legend6SubGrad03}{\scriptsize SubGD: $\zeta^k=10^{-3}/k^{0.3}$}
\psfrag{Objective}[cc][Bt]{Objective function}
\psfrag{Itr}[cc][Bl]{Iteration}
\includegraphics[width=0.24\textwidth]{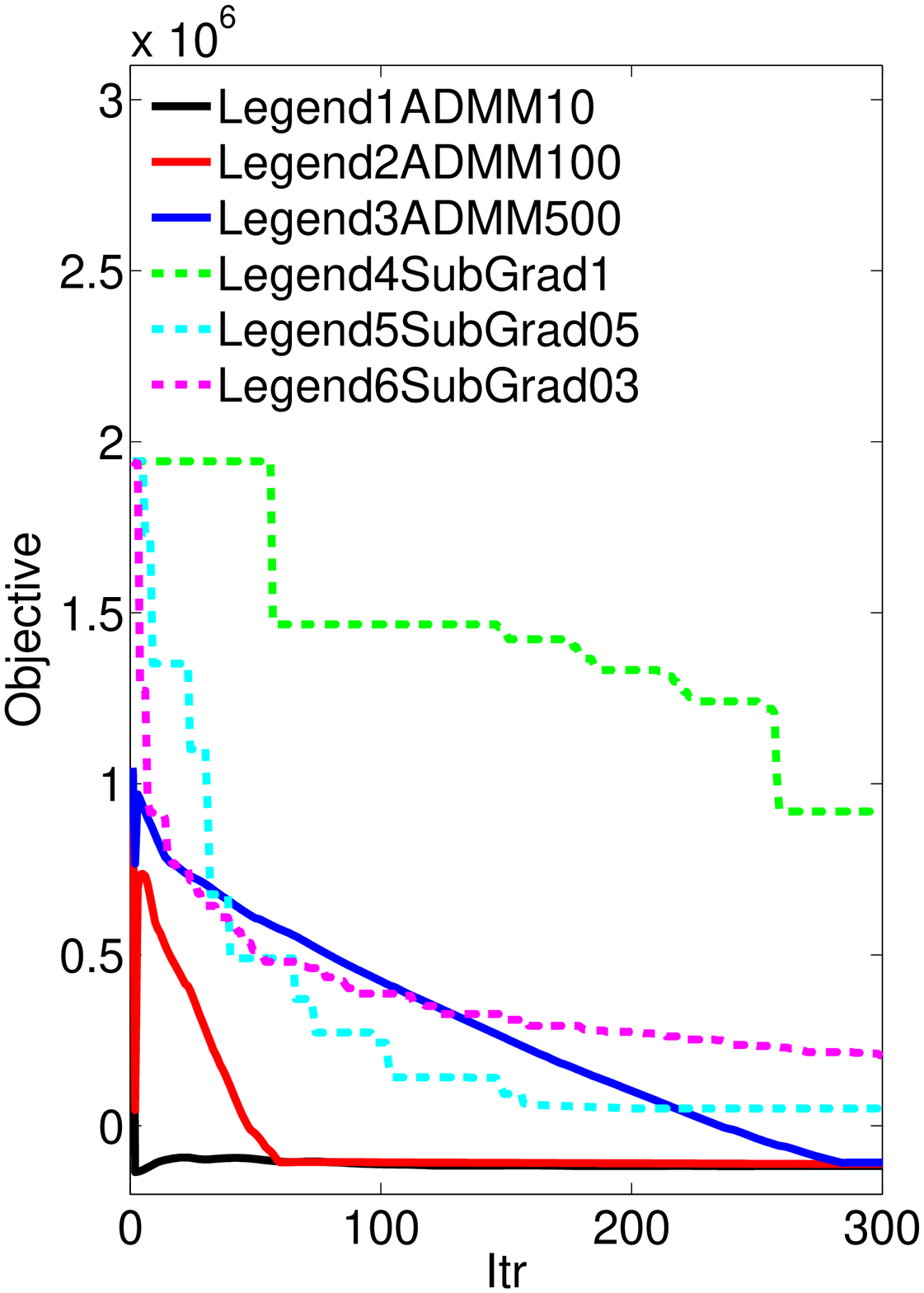}}
\centering
\subfigure[ADMM residual convergence]{ \label{fig:ADMM_feas}
\centering
\psfrag{Legend1ADMM10}{\scriptsize ADMM: $\rho=10$}
\psfrag{Legend2ADMM100}{\scriptsize ADMM: $\rho=100$}
\psfrag{Legend3ADMM500}{\scriptsize ADMM: $\rho=500$}
\psfrag{Legend4SubGrad1}{\scriptsize SubGD: $\zeta^k=10^{-3}/k$}
\psfrag{Legend5SubGrad05}{\scriptsize SubGD: $\zeta^k=10^{-3}/k^{0.5}$}
\psfrag{Legend6SubGrad03}{\scriptsize SubGD: $\zeta^k=10^{-3}/k^{0.3}$}
\psfrag{Residual}[cc][Bt]{Primal residual}
\psfrag{Itr}[cc][Bl]{Iteration}
\includegraphics[width=0.24\textwidth]{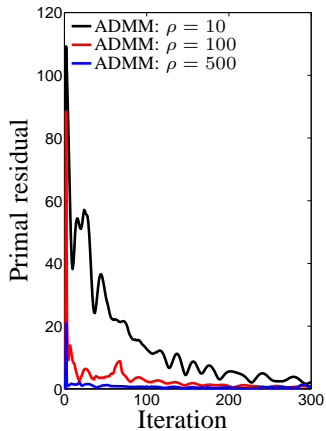}}
\caption{Convergence of ADMM and centralized subgradient method. Here $\zeta^k$ is the step size of the subgradient algorithm at the $k$th iteration.}\label{fig:6}
\end{figure}

The convergence of the fully distributed ADMM is shown in Figure~\ref{fig:6}. As a comparison, we also plot the convergence of the projected subgradient method (SubGD). Figure~\ref{fig:ADMM_subgrad_conv} shows the convergence of the objective value of the online program at a time period for both algorithms with different algorithmic parameter choices, while Figure~\ref{fig:ADMM_feas} depicts the convergence of the norm of the primal residual for the ADMM algorithm. In terms of the objective value, we observe that the convergence of ADMM is usually much faster compared to SubGD. In fact, in all our examples, SubGD does not converge  after thousands of iterations with the tolerance being $1\times 10^{-4}$. Comparing the performance of ADMM with different parameter $\rho$'s, we note that smaller $\rho$ leads to faster convergence in terms of the objective value but slower convergence of the primal residual. Thus in practice, selecting a $\rho$ that properly trades off these two effects is necessary. 

\subsection{Scalability}\label{sec:nt:scale}
\begin{figure}[htbp]
\
\centering
\psfrag{ADMMWithSingleProcessor}{\scriptsize ADMM with single processor}
\psfrag{ADMMWithMultiProcessors}{\scriptsize Fully distributed ADMM}
\psfrag{RunTimeSeconds}[cc][Bt]{Running time (second)}
\psfrag{NumberOfNodes}{Number of nodes}
\psfrag{15000}{ \scriptsize 15000}
\psfrag{10000}{ \scriptsize 10000}
\psfrag{5000}{ \scriptsize 5000}
\psfrag{0}[cc][Bt]{ \scriptsize 0}
\psfrag{5}{ \scriptsize 5}
\psfrag{10}{ \scriptsize 10}
\psfrag{15}{ \scriptsize 15}
\includegraphics[width=0.4\textwidth]{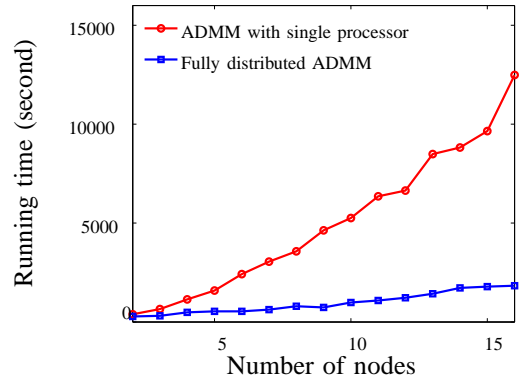}
\caption{Running time of distributed and centralized ADMM}\label{fig:ADMM_scale}
\end{figure}

In this subsection, we give a preliminary account for the scalability of the distributed implementation using Matlab Distributed Computing Toolbox. We consider star networks discussed in Subsection~\ref{sec:nt:star} with the number of buses increasing from $2$ to $16$. We associate a processor to each of the buses, and run the distributed ADMM using $2$-$16$ processors. The running time, together with the running time of solving the online programs using the centralized ADMM algorithm, is shown in Figure~\ref{fig:ADMM_scale}. We note that while in both scenarios, the running time increases approximately linearly with the number of buses, the rate of linear increase for distributed ADMM is significantly smaller. Loading the data for problem specification and communication overheads may have contributed to the linear running time increase for the distributed ADMM.

\section{Conclusion and Future Directions}
This paper formulates the storage network operation problem as a stochastic control problem. An online algorithm is proposed to solve the problem efficiently. The performance of the algorithm is analyzed and a sub-optimality bound is derived. The online programs are then solved in a decentralized fashion with only local computation and neighborhood communication with task-based ADMM iterations. Combining these elements, we obtain an efficient task-based distributed online control strategy for operating distributed storage systems with a guaranteed performance.

Many future directions are of interest for generalizing our results. (i) This paper focuses on the real power; incorporating the reactive power and a full AC power flow model may be an important step towards a successful implementation in large-scale practical systems. As the online optimization for each step becomes an AC optimal power flow (OPF) problem, recent work on the convexification of such problems \cite{slowOPF1} \cite{slowOPF2}, and the distributed solution of the convexified program \cite{6502290} may be combined with the approach proposed in this paper. (ii) Our decentralized solution is based on the classical two block ADMM which has superior convergence properties compared to other popular methods for distributed optimization such as primal or dual decomposition. Similar methods have been tested in much larger networks for deterministic energy control problems \cite{kraning2013dynamic}. However, the fact that such an ADMM algorithm requires a two-block partition (corresponding to the node variables $x$ and edge variables $z$ in Section~\ref{sec:admm:n+e}) leads to the inconvenience that local copies of variables controlled by the neighbors must be created. Multi-block variants of ADMM may eliminate such need. However, the convergence is not guaranteed or requires additional assumptions \cite{admm3b}, \cite{hong2012linear}, \cite{han2012note}. Validating these assumptions for specific storage control problem instances may lead to simpler algorithm which has similar convergence properties. (iii) Utilizing the sub-optimality bounds to assess the limit of the storage system for the purpose of storage valuation and system design may also be of interest.


\bibliography{jqin}

\begin{thebibliography}{10}
\providecommand{\url}[1]{#1}
\csname url@samestyle\endcsname
\providecommand{\newblock}{\relax}
\providecommand{\bibinfo}[2]{#2}
\providecommand{\BIBentrySTDinterwordspacing}{\spaceskip=0pt\relax}
\providecommand{\BIBentryALTinterwordstretchfactor}{4}
\providecommand{\BIBentryALTinterwordspacing}{\spaceskip=\fontdimen2\font plus
\BIBentryALTinterwordstretchfactor\fontdimen3\font minus
  \fontdimen4\font\relax}
\providecommand{\BIBforeignlanguage}[2]{{%
\expandafter\ifx\csname l@#1\endcsname\relax
\typeout{** WARNING: IEEEtran.bst: No hyphenation pattern has been}%
\typeout{** loaded for the language `#1'. Using the pattern for}%
\typeout{** the default language instead.}%
\else
\language=\csname l@#1\endcsname
\fi
#2}}
\providecommand{\BIBdecl}{\relax}
\BIBdecl

\bibitem{NRELWest2010}
\BIBentryALTinterwordspacing
{National Renewable Energy Laboratory}. (2010) {Western Wind and Solar
  Integration Study}. [Online]. Available:
  \url{http://www.nrel.gov/wind/systemsintegration/wwsis.html}
\BIBentrySTDinterwordspacing

\bibitem{Denholm2010}
\BIBentryALTinterwordspacing
------. (2010) {The Role of Energy Storage with Renewable Electricity
  Generation}. [Online]. Available:
  \url{http://www.nrel.gov/wind/pdfs/47187.pdf}
\BIBentrySTDinterwordspacing

\bibitem{lindley2010naturenews}
D.~Lindley, ``{Smart Grids: The Energy Storage Problem},'' \emph{Nature}, vol.
  463, no. 7277, p.~18, 2010.

\bibitem{thermalStor1993}
B.~Daryanian and R.~E. Bohn, ``{Sizing of Electric Thermal Storage under Real
  Time Pricing},'' \emph{{IEEE} Transactions on Power Systems}, vol.~8, no.~1,
  pp. 35--43, 1993.

\bibitem{ObRACC2013}
G.~O'Brien and R.~Rajagopal, ``{A Method for Automatically Scheduling Notified
  Deferrable Loads},'' in \emph{Proc. of American Control Conference (ACC)},
  2013, pp. 5080--5085.

\bibitem{HaoSanandajiPoollaVincent2013}
H.~Hao, B.~M. Sanandaji, K.~Poolla, and T.~L. Vincent, ``{Aggregate Flexibility
  of Thermostatically Controlled Loads},'' \emph{IEEE Transactions on Power
  Systems}, submitted.

\bibitem{QRsimpleStorPes2012}
J.~Qin, R.~Sevlian, D.~Varodayan, and R.~Rajagopal, ``{Optimal Electric Energy
  Storage Operation},'' in \emph{Proc. of IEEE Power and Energy Society General
  Meeting}, 2012, pp. 1--6.

\bibitem{MITrampStor}
\BIBentryALTinterwordspacing
A.~{Faghih}, M.~{Roozbehani}, and M.~A. {Dahleh}, ``{On the Economic Value and
  Price-Responsiveness of Ramp-Constrained Storage},'' \emph{{ArXiv} e-prints},
  2012. [Online]. Available: \url{arxiv.org/abs/1211.1696}
\BIBentrySTDinterwordspacing

\bibitem{SuEGTPS}
H.~I. Su and A.~El~Gamal, ``{Modeling and Analysis of the Role of Energy
  Storage for Renewable Integration: Power Balancing},'' \emph{IEEE
  Transactions on Power Systems}, vol.~28, no.~4, pp. 4109--4117, 2013.

\bibitem{RLDSACC}
J.~Qin, H.~I. Su, and R.~Rajagopal, ``{Storage in Risk Limiting Dispatch:
  Control and Approximation},'' in \emph{Proc. of American Control Conference
  (ACC)}, 2013, pp. 4202--4208.

\bibitem{BitarRACC_colocated}
E.~Bitar, R.~Rajagopal, P.~Khargonekar, and K.~Poolla, ``{The Role of
  Co-Located Storage for Wind Power Producers in Conventional Electricity
  Markets},'' in \emph{Proc. of American Control Conference (ACC)}, 2011, pp.
  3886--3891.

\bibitem{Powell}
\BIBentryALTinterwordspacing
J.~H. Kim and W.~B. Powell, ``{Optimal Energy Commitments with Storage and
  Intermittent Supply},'' \emph{Operations Research}, vol.~59, no.~6, pp.
  1347--1360, 2011. [Online]. Available:
  \url{http://pubsonline.informs.org/doi/abs/10.1287/opre.1110.0971}
\BIBentrySTDinterwordspacing

\bibitem{IBMload}
\BIBentryALTinterwordspacing
P.~M. {van de Ven}, N.~{Hegde}, L.~{Massoulie}, and T.~{Salonidis}, ``{Optimal
  Control of End-User Energy Storage},'' \emph{{ArXiv} e-prints}, 2012.
  [Online]. Available: \url{arxiv.org/abs/1203.1891}
\BIBentrySTDinterwordspacing

\bibitem{DataCenter}
R.~Urgaonkar, B.~Urgaonkar, M.~J. Neely, and A.~Sivasubramaniam, ``{Optimal
  Power Cost Management Using Stored Energy in Data Centers},'' in \emph{Proc.
  of the {ACM} SIGMETRICS Joint International Conference on Measurement and
  Modeling of Computer Systems}, ser. SIGMETRICS '11, 2011, pp. 221--232.

\bibitem{StorDRLongbo}
L.~Huang, J.~Walrand, and K.~Ramchandran, ``{Optimal Demand Response with
  Energy Storage Management},'' in \emph{Proc. of IEEE Third International
  Conference on Smart Grid Communications (SmartGridComm)}, 2012, pp. 61--66.

\bibitem{XieEtAlWindStorMPC}
L.~Xie, Y.~Gu, A.~Eskandari, and M.~Ehsani, ``{Fast MPC-Based Coordination of
  Wind Power and Battery Energy Storage Systems},'' \emph{Journal of Energy
  Engineering}, vol. 138, no.~2, pp. 43--53, 2012.

\bibitem{NRELStorValue2013}
\BIBentryALTinterwordspacing
{National Renewable Energy Laboratory}. (2013) {The Value of Energy Storage for
  Grid Applications}. [Online]. Available:
  \url{http://www.nrel.gov/docs/fy13osti/58465.pdf}
\BIBentrySTDinterwordspacing

\bibitem{OMGarxiv}
\BIBentryALTinterwordspacing
J.~{Qin}, Y.~{Chow}, J.~{Yang}, and R.~{Rajagopal}, ``{Online Modified Greedy
  Algorithm for Storage Control under Uncertainty},'' \emph{{ArXiv} e-prints},
  2014. [Online]. Available: \url{arxiv.org/abs/1405.7789}
\BIBentrySTDinterwordspacing

\bibitem{QCYR:acm}
J.~Qin, Y.~Chow, J.~Yang, and R.~Rajagopal, ``{Modeling and Online Control of
  Generalized Energy Storage Networks},'' in \emph{Proc. of the 5th
  International Conference on Future Energy Systems (ACM e-Energy '14)}.\hskip
  1em plus 0.5em minus 0.4em\relax ACM, June 2014.

\bibitem{slowOPF1}
\BIBentryALTinterwordspacing
S.~H. {Low}, ``{Convex Relaxation of Optimal Power Flow, Part I: Formulations
  and Equivalence},'' \emph{ArXiv e-prints}, May 2014. [Online]. Available:
  \url{arxiv.org/abs/1405.0766}
\BIBentrySTDinterwordspacing

\bibitem{slowOPF2}
\BIBentryALTinterwordspacing
------, ``{Convex Relaxation of Optimal Power Flow, Part II: Exactness},''
  \emph{ArXiv e-prints}, May 2014. [Online]. Available:
  \url{arxiv.org/abs/1405.0814}
\BIBentrySTDinterwordspacing

\bibitem{QCYR:report}
\BIBentryALTinterwordspacing
J.~Qin, Y.~Chow, J.~Yang, and R.~Rajagopal, ``{Control of Generalized Energy
  Storage Networks},'' \emph{{Stanford S3L Report}}, 2014, available at
  \url{http://www.stanford.edu/~jqin/pdf/QCYR2014report.pdf}. [Online].
  Available: \url{http://www.stanford.edu/~jqin/pdf/QCYR_2013_report.pdf}
\BIBentrySTDinterwordspacing

\bibitem{grip:cdc2011}
D.~Bakken, A.~Bose, K.~Chandy, P.~Khargonekar, A.~Kuh, S.~Low, A.~Von~Meier,
  K.~Poolla, P.~Varaiya, and F.~Wu, ``Grip - grids with intelligent periphery:
  Control architectures for grid2050,'' in \emph{Smart Grid Communications
  (SmartGridComm), 2011 IEEE International Conference on}, Oct 2011, pp. 7--12.

\bibitem{BoydADMM}
S.~Boyd, N.~Parikh, E.~Chu, B.~Peleato, and J.~Eckstein, ``{Distributed
  Optimization and Statistical Learning Via the Alternating Direction Method of
  Multipliers},'' \emph{Foundations and Trends{\textregistered} in Machine
  Learning}, vol.~3, no.~1, pp. 1--122, 2011.

\bibitem{wei20131}
\BIBentryALTinterwordspacing
E.~Wei and A.~Ozdaglar, ``{On the $O(1/k)$ Convergence of Asynchronous
  Distributed Alternating Direction Method of Multipliers},'' \emph{{ArXiv}
  e-prints}, 2013. [Online]. Available: \url{arxiv.org/abs/1307.8254}
\BIBentrySTDinterwordspacing

\bibitem{IEEEtest}
\BIBentryALTinterwordspacing
Power systems test case archive. [Online]. Available:
  \url{www.ee.washington.edu/research/pstca/}
\BIBentrySTDinterwordspacing

\bibitem{NREL2010}
{National Renewable Energy Laboratory}, ``Eastern wind integration and
  transmission study,'' Tech. Rep., 2010.

\bibitem{6502290}
E.~Dall'Anese, H.~Zhu, and G.~Giannakis, ``{Distributed Optimal Power Flow for
  Smart Microgrids},'' \emph{IEEE Transactions on Smart Grid}, vol.~4, no.~3,
  pp. 1464--1475, Sept 2013.

\bibitem{kraning2013dynamic}
M.~Kraning, E.~Chu, J.~Lavaei, and S.~Boyd, ``{Dynamic Network Energy
  Management Via Proximal Message Passing},'' \emph{Foundations and Trends in
  Optimization}, vol.~1, no.~2, pp. 1--54, 2013.

\bibitem{admm3b}
\BIBentryALTinterwordspacing
C.~Chen, B.~S. He, Y.~Ye, and X.~Yuan, ``{The Direct Extension of ADMM for
  Multi-Block Convex Minimzation Problems Is Not Necessarily Convergent},''
  \emph{Mathematical Programming}, to appear. [Online]. Available:
  \url{www.optimization-online.org/DB_FILE/2013/09/4059.pdf}
\BIBentrySTDinterwordspacing

\bibitem{hong2012linear}
\BIBentryALTinterwordspacing
M.~Hong and Z.-Q. Luo, ``{On the Linear Convergence of the Alternating
  Direction Method of Multipliers},'' \emph{{arXiv} e-print}, 2012. [Online].
  Available: \url{arxiv.org/abs/1208.3922}
\BIBentrySTDinterwordspacing

\bibitem{han2012note}
D.~Han and X.~Yuan, ``{A Note on the Alternating Direction Method of
  Multipliers},'' \emph{Journal of Optimization Theory and Applications}, vol.
  155, no.~1, pp. 227--238, 2012.

\bibitem{NeelyBook}
M.~J. Neely, ``{Stochastic Network Optimization with Application to
  Communication and Queueing Systems},'' \emph{Synthesis Lectures on
  Communication Networks}, vol.~3, no.~1, pp. 1--211, 2010.

\bibitem{admmRate}
B.~He and X.~Yuan, ``{On the $O(1/n)$ Convergence Rate of the Douglas--Rachford
  Alternating Direction Method},'' \emph{SIAM Journal on Numerical Analysis},
  vol.~50, no.~2, pp. 700--709, 2012.

\bibitem{2012arXiv1208.3922H}
\BIBentryALTinterwordspacing
M.~{Hong} and Z.-Q. {Luo}, ``{On the Linear Convergence of the Alternating
  Direction Method of Multipliers},'' \emph{ArXiv e-prints}, Aug. 2012.
  [Online]. Available: \url{arxiv.org/abs/1208.3922}
\BIBentrySTDinterwordspacing

\end{thebibliography}
\bibliographystyle{IEEEtran}

\appendices

\section{Definitions and Expressions for Section~\ref{sec:omgnet}}\label{app:def4omg}
Here we provide the definition and expressions for $\Dl \g_i$, $\Du \g_i$, $\ksmin_i$, $\ksmax_i$ and $\Wmax_i$. We start by defining $\Dl \g_i$ and $\Du \g_i$ for each $i \in [n]$. 
\begin{definition}
Let $y_i\defeq (\f,\d_i, \p_i)$. For function $\phi_i(u_i, y_i) \defeq \g_i(\u_i-\d_i + \sum_{j=1}^n  \f_{ij}, \p_i)$  that is convex (but not necessarily differentiable) in $\u_i$, a real number $\alpha_i$ is called a (partial) subgradient of $\phi_i$ with respect to argument $\u_i$ at given $(\u_i, y_i)$ if $\phi_i(u_i',y_i) \ge \phi_i(u_i, y_i) + \alpha_i (u_i'-u_i)$ for all $u_i' \in [\umin_i, \umax_i]$. The set of all subgradients at $(\u_i, y_i)$, denoted by $\partial_{u_i} \phi_i(u_i,y_i)$, is called the (partial) subdifferential of $\phi_i(u_i, y_i)$ with respect to $\u_i$ at $(\u_i,y_i)$. Denote $\mathcal{U}_i \defeq [\umin_i,\umax_i]$, $\mathcal{Y}_i \defeq \Fcal \times [\dmin_i, \dmax_i] \times [\pmin_i, \pmax_i]$ where $\Fcal=\big\{\f: -\Fmax_{ij}\le\f_{ij}  \le \Fmax_{ij}, \,\forall i,j \in [n]\big\}$. Define the set
\[
D \g_i \defeq \bigcup_{(u_i,y_i) \in \times\mathcal{U}_i \times \mathcal{Y}_i} \partial_{u_i} \phi_i(\u_i, y_i),
\]
and let real numbers $\Dl g_i$ and $\Du g_i$ be defined such that
\begin{equation}
\Dl\g_i \le \inf D\g_i \le \sup D\g_i \le \Du\g_i.
\end{equation}
That is, $\Dl \g_i$ and $\Du \g_i$ are a lower bound and an upper bound of the sub-gradient of $\phi_i$ over  its (compact) domain, respectively.
\end{definition}
More details and examples regarding how to calculate $\Dl \g_i$ and $\Du \g_i$ can be found in our previous work \cite{OMGarxiv}. The bounds for the algorithmic parameters are
\begin{equation}\label{eq:ineq_1}
\ksmin_i \defeq \frac{1}{\la_i} \left(-\W_i \Dl \g_i + \umax_i - \bmax_i\right),
\end{equation}
\begin{equation}\label{eq:ineq_2}
\ksmax_i \defeq \frac{1}{\la_i} \left(-\W_i \Du \g_i - \bmin_i+ \umin_i\right),\end{equation}
and
\begin{equation}\label{eq:W_max}
\Wmax_i \defeq \frac{(\bmax_i - \bmin_i) -(\umax_i - \umin_i)}{\Du \g_i - \Dl \g_i}  .
\end{equation}
\section{Proof of Theorem \ref{thm:perf_lyap}}\label{app:lyap_pf}
Similar to the analysis in \cite{OMGarxiv} for the single bus storage case, we will prove Theorem \ref{thm:perf_lyap} via the following steps:
\begin{enumerate}
\item Reformulate problem~\eqref{p1_inf} and link it to the sequence of OMG online optimizations \eqref{OMGcen}.
\item Prove that the control policy obtained from OMG is feasible to problem~\eqref{p1_inf}.
\item Derive the performance bound in Theorem \ref{thm:perf_lyap}.
\end{enumerate}

First, we proceed by reformulating problem~\eqref{p1_inf}.
For $i=1,\ldots,n$, define
\[
\bar \u_i \defeq \lim_{T \to \infty} \frac{1}{T} \expec\left[ \sum_{t=1}^T \u_i(t)\right],\,\,\bar \b_i \defeq \lim_{T \to \infty} \frac{1}{T} \expec\left[ \sum_{t=1}^T \b_i(t)\right].
\]
Note that for $\b_i(1)\in[\bmin_i,\bmax_i]$,
\[
\bar \u_i=\lim_{T \to \infty}\frac{1}{T}\expec\left[\sum_{t=1}^T\b_i(t+1)-\la_i\b_i(t)\right]=(1-\la_i)\bar\b_i.
\]
As $\b_i(t)\in[\bmin_i,\bmax_i]$ for all $t\geq 0$, the above expression implies
\[
(1-\la_i)\bmin_i\leq \bar \u_i\leq (1-\la_i)\bmax_i.
\]
Problem~\eqref{p1_inf} can be equivalently written as follows
\begin{subequations}\label{prob:P1}
\begin{align}
\opttag{P1:} \minimize &  \lim_{T\to \infty}\frac{1}{T}\expec \sum_{t=1}^T \sum_{i=1}^n \g_{i}(\r_i(t); \p_i(t)) \\
\st & \delta_i(t) +\r_i(t)  = \u_i(t) + \sum_{j=1}^n  \f_{ij}(t) , \label{P1:balance}\\
 &\b_i(t+1) = \la_i \b_i(t) + \u_i(t), \label{P1:dynamics}\\
        & \bmin_i - \la_i \b_i(t) \le \u_i(t) \le \bmax_i - \la_i \b_i(t), \label{P1:bbounds-u}\\
        & \umin_i \le \u_i(t) \le \umax_i, \label{P1:ubounds}\\
        &(1-\la_i)\bmin\leq \bar \u_i\leq (1-\la_i)\bmax_i \label{P1:ubar}\\
        &B_{ij} (\theta_i(t) -\theta_j(t)) = \f_{ij}(t) ,\label{P1:theta}\\
& -\Fmax_{ij} \le \f_{ij}(t) \le \Fmax_{ij} \label{P1:F},
\end{align}
\end{subequations}
where bounds on $\b_i(t)$ are replaced by \eqref{P1:bbounds-u}, and \eqref{P1:ubar} is added without loss of optimality.

Here we use  $\J_\mathrm{P1}(\u,\r,\theta,\f)$ to denote the objective value of {\bf P1} with operation sequence $(\u,\r,\theta,\f)$ (as an abbreviation of $\{\u(t),\r(t),\theta(t),\f(t): t\ge 1\}$), $\Lambda^\star(\mathrm{\bf P1})=(\u^\star(\mathrm{\bf P1}),\r^\star(\mathrm{\bf P1}),\theta^\star(\mathrm{\bf P1}),\f^\star(\mathrm{\bf P1}))$ to denote the optimal control sequence for {\bf P1}, $\J^\star_\mathrm{P1} \defeq  \J_\mathrm{P1}(\Lambda^\star(\mathrm{\bf P1}))$, and we define similar quantities for {\bf P2}. Here {\bf P2} is an auxilliary problem we construct to bridge the infinite horizon storage control problem {\bf P1} to online optimization problems~\eqref{OMGcen}. It has the following form
\begin{subequations}\label{prob:P2}
\begin{align}
\opttag{P2:} \minimize  &  \lim_{T \to \infty} \frac{1}{T}  \expec \sum_{t=1}^T \sum_{i=1}^n \g_{i}(\r_i(t); \p_i(t))\\
\st  
       &  \delta_i(t) +\r_i(t)  = \u_i(t) + \sum_{j=1}^n  \f_{ij}(t) , \label{P2:balance}\\
        & \umin_i \le \u_i(t) \le \umax_i, \label{P2:ubounds}\\
        &(1-\la_i)\bmin_i\leq \bar \u_i\leq (1-\la_i)\bmax_i \label{P2:ubar}\\
        &B_{ij} (\theta_i(t) -\theta_j(t)) = \f_{ij}(t) ,\label{P2:theta}\\
& -\Fmax_{ij} \le \f_{ij}(t) \le \Fmax_{ij} \label{P2:F}.
\end{align}
\end{subequations}
Notice that it has the same objective as \Pk{1}, and evidently it is a relaxation of \Pk{1}. This implies that $\u^\star(\Pk{2})$ may not be feasible for \Pk{1}, and
\begin{equation}
\Jk{2}^\star = \Jk{1}(\Lambda^\star(\Pk{2})) \le \Jk{1}^\star.
\end{equation}
The reason for the removal of state-dependent constraints \eqref{P1:bbounds-u} (and hence \eqref{P1:dynamics} as the sequence $\{\b(t): t\ge 1\}$ becomes irrelevant to the optimization of $\{\u(t): t\ge 1\}$) in \Pk{2} is that the state-independent problem \Pk{2} has easy-to-characterize optimal stationary control policies. In particular, from the theory of stochastic network optimization \cite{NeelyBook}, the following result holds.
\begin{lemma}
[Stationary Disturbance-Only Policies]\label{lem:stat_network}
Under Assumption~\ref{ass:1}
there exists a stationary disturbance-only policy $\Lambda^{\mathrm{stat}}(t)=(\ustatn(t),\rstatn(t),\thetastatn(t),\fstatn(t))$ satisfying the constraints in \Pk{2} and providing the following guarantees $\forall t$:
\begin{align}
  & (1-\la_i)\bmin_i\leq \expec[\ustat_i(t) ] \leq (1-\la_i)\bmax_i, \forall i\in [i]\nn \\
  & \expec\left[ \sum_{i=1}^n \g_{i}(\r_i(t); \p_i(t))\middle| \Lambda^{\mathrm{stat}}(t) \right]  = \Jk{2}^\star,\nn\end{align}
where the expectation is taken over the randomization of $\d_i(t)$, $\p_i(t)$, and possibly $\Lambda^{\mathrm{stat}}(t)$ in case the policy is randomized.
\end{lemma}

Recall the online optimization solved by OMG:
\begin{subequations}\label{OMGP3}
\begin{align}
\opttag{P3:} \minimize  &  \sum_{i=1}^n (\la_i/\W_i) (\b_i+\ks_i) \u_i + \g_{i}(\r_i; \p_i)\\
\st  
        & \umin_i \le u_i \le \umax_i,\\
        & \delta_i +\r_i  = \u_i + \sum_{j=1}^n  \f_{ij},\\
        &B_{ji} (\theta_j -\theta_i) = \f_{ji} ,\\
        & -\Fmax_{ji} \le \f_{ji} \le \Fmax_{ji}.
\end{align}
\end{subequations}
We use $\Lambda^{\mathrm{ol}}(t)=(\uhatn(t),\rhatn(t),\thetahatn(t),\fhatn(t))$ to denote the solution of \Pk{3} at time step $t$, 
$\Lambda^\star(\Pk{3})=(\u^\star(\Pk{3}),\r^\star(\Pk{3}),\theta^\star(\Pk{3}),\fn^\star(\Pk{3})) $
to denote the sequence $\{\Lambda^{\mathrm{ol}}(t): t\ge 1\}$, $\J_{\Pk{3},t}(\Lambda(t))$ to denote the objective function of \Pk{3} at time period $t$ using policy $\Lambda(t)$, and $\J_{\Pk{3},t}^\star$ to denote the corresponding optimal cost.

Now, we turn to the feasibility analysis of $\Lambda^\star(\Pk{3})$ with respect to \Pk{1}. Following assumption holds for any storage system that is \emph{controllable}. 
\begin{assumption}[Feasibility and Controllability]\label{assume:feas}
Each storage $i\in[n]$ is feasible and controllable: 
\begin{itemize} 
\item (feasibility) starting from any feasible storage level, there exists a feasible storage operation such that the storage level in the next time period is feasible, \ie, $\la_i \bmin_i + \umax_i \ge \bmin_i$ and $\la_i \bmax_i + \umin_i \le \bmax_i$.
\item (controllability) starting from any feasible storage level, 
there exists a sequence of feasible storage operations to reach any feasible storage level in a finite number of steps, \ie, $\la_i \bmax_i + \umax_i \ge \bmax_i$ and $\la_i \bmin_i + \umin_i \le \bmin_i$.
\end{itemize}
%
%
%
\end{assumption}
In order to prove that the solution of \Pk{3} is feasible to \Pk{1}, we have the following technical lemma. 
\begin{lemma}\label{coro:tech_res}
At each time period $t$, the optimal storage operation of \Pk{3} at node $i$, $\uhat_i(t)$, for $i=1,\ldots,n$, satisfies
\begin{enumerate}
\item $\uhat_i(t) = \umin_i$ whenever $\la_i \bs_i(t) \ge - \W_i \Dl\g_i $,
\item $\uhat_i(t) = \umax_i$ whenever $\la_i \bs_i(t) \le - \W_i \Du\g_i $,
\end{enumerate}
where 
\[
\bs_i(t) = \b_i(t) + \ks_i.
\]
\end{lemma}
\begin{IEEEproof}
The proof follows from similar arguments used to prove Lemma 3 of \cite{OMGarxiv}. Details are omitted for brevity. 
\end{IEEEproof}

We are ready to prove the feasibility of the control sequence generated by the algorithm. 
\begin{IEEEproof}[Proof of Theorem~\ref{thm:perf_lyap}, feasibility]
For any $i=1,\ldots,n$, we first validate that the intervals of $\ks_i$ and $\W_i$ are non-empty. By {\bf A3} of Assumption \ref{ass:1}, one concludes $\Wmax_i> 0$, thus it remains to show $\ksmax_i \ge \ksmin_i$. 
Based on \eqref{eq:W_max}, $W_i\geq0$, and $\Du\g_i\geq \Dl\g_i$, one obtains
\[
\W_i(\Du \g_i - \Dl \g_i)\leq [(\bmax_i - \bmin_i) - (\umax_i-\umin_i)].
\]
Re-arranging terms results in
\[
-\W_i \Dl g_i + \umax_i - \bmax_i \le - \W_i \Du g_i - \bmin_i + \umin_i,
\]
which further implies $\ksmax_i\geq \ksmin_i$.

We proceed to show that
\begin{equation}\label{exp:MI}
\bmin_i\leq  \b_i(t)\leq \bmax_i,
\end{equation}
for $t= 1,2,\dots$ and any $i \in [n]$, when $\Lambda^\star(\Pk{3})$ is implemented. 
The base case holds by assumption.
Let the inductive hypothesis be that (\ref{exp:MI}) holds at time $t$.
The storage level at $t+1$ is  then
$\b_i(t+1)=\la_i \b_i(t)+\uhat_i(t). $
We show  (\ref{exp:MI}) holds at $t+1$ by considering the following three cases.

\noindent{\bf Case 1.}  $- \W_i \Dl\g_i\leq\la_i\bs_i(t)\leq\la_i ( \bmax_i+\ks_i)$. \\
First, it is easy to verify that the above interval for $\la_i \bs_i(t)$ is non-empty using \eqref{eq:ineq_1} and $\ks_i \ge \ksmin_i$. Next, based on Lemma \ref{coro:tech_res}, one obtains  $\uhat_i(t)=\umin_i\le 0$ in this case. Therefore
\[
\b_i(t+1) = \la_i \b_i(t) + \umin_i \le  \la_i \bmax_i + \umin_i \le \bmax_i,
\]
where the last inequality follows from Assumption \ref{assume:feas}. On the other hand, 
\begin{align*}
\b_i(t+1)& = \la_i \b_i(t) + \umin_i \ge -\W_i\Dl \g_i  - \la_i \ks_i + \umin_i \\
\ge & -\W_i\Dl \g_i  - \la_i \ksmax_i + \umin_i 
\ge    \bmin_i, 
\end{align*}
where the third inequality used $\Du \g_i \ge \Dl \g_i$. 

\noindent{\bf Case 2.}
 $\la_i( \bmin_i+\ks_i) \le  \la_i\bs_i(t)  \le - \W_i \Du\g_i$.  \\
The above interval for $\la_i \bs_i(t)$ is non-empty by \eqref{eq:ineq_2} and $\ks_i \le \ksmax_i$. Lemma \ref{coro:tech_res} implies  $\uhat_i(t)=\umax_i \ge 0$ in this case. 
Therefore, again using Assumption \ref{assume:feas},
\[
\b_i(t+1) = \la_i \b_i(t) + \umax_i \ge \la_i \bmin_i + \umax_i \ge \bmin_i.
\]
On the other hand, 
\begin{align*}
\!\!\b_i(t+1) &= \la_i \b_i(t) + \umax_i \le -\W_i\Du \g_i - \la_i \ks_i + \umax_i \\
\le &  -\W_i\Du \g_i - \la_i \ksmin_i + \umax_i 
\le \bmax_i,
\end{align*}
where the third inequality used $\Du \g_i \ge \Dl \g_i$. 

\noindent{\bf Case 3.} $-\W_i \Du\g_i<\la_i\bs_i(t)<- \W_i \Dl\g_i$.  \\
By  $\umin_i \le \uhat_i(t)\le\umax_i$, one obtains
\begin{align*}
\b_i(t+1) &=  \la_i \b_i(t) + \uhat_i(t) \le \la_i \b_i(t) + \umax_i\\
< & - \W_i \Dl\g_i - \la_i \ks_i + \umax_i\\
\le &  - \W_i \Dl\g_i - \la_i \ksmin_i + \umax_i
\le \bmax_i.
\end{align*}
On the other hand,
\begin{align*}
\b_i(t+1) &=  \la_i \b_i(t) + \uhat_i(t) \ge \la_i \b_i(t) + \umin_i\\
> & - \W_i \Du\g_i - \la_i \ks_i + \umin_i\\
\ge &  - \W_i \Du\g_i - \la_i \ksmax_i + \umax_i
\ge    \bmin_i.
\end{align*}

Combining these three cases, and by mathematical induction, we conclude \eqref{exp:MI} holds for all $t = 1,2,\dots$. 
\end{IEEEproof}

It remains to show that the sub-optimality bounds claimed in Theorem~\ref{thm:perf_lyap} indeed hold. 
\begin{IEEEproof}[Proof of Theorem~\ref{thm:perf_lyap}, performance]
Consider a quadratic Lyapunov function $\L_i(\b_i) = \b_i^2/2$.
Let the corresponding Lyapunov drift be
\[
\Ls_i( \b_i(t)) = \expec\left[ \L_i(\b_i(t+1)) -  \L_i(\b_i(t)) \vert  \b_i(t) \right].
\]
Recall that
$
\bs_i(t+1) = \b_i(t+1) + \ks_i = \la \bs_i(t)   + \u_i(t) + (1-\la_i)\ks_i,
$
and so
\begin{align}
&\quad\,\,\Ls_i( \bs_i(t))\nn\\
& = \expec\big[ (1/2)(\u_i(t) + (1-\la_i)\ks_i)^2  -  (1/2) (1-\la_i^2)\bs_i(t)^2 \nn\\
                &\quad\quad\quad + \la_i \bs_i(t) \u_i(t) + \la_i (1-\la_i) \bs_i(t) \ks_i \vert  \bs_i(t) \big]\nn\\
& \le \Mone_i(\ks_i)    -  (1/2) (1-\la_i^2)\bs_i(t)^2\nn \\
& \quad+\expec\big[ \la_i \bs_i(t) \u_i(t) + \la_i (1-\la_i) \bs_i(t) \ks_i \vert  \bs_i(t) \big]\nn\\
& \le\Mone_i(\ks_i) + \expec\left[ \la_i \bs_i(t) (\u_i(t)+(1-\la_i)\ks_i)\vert \bs_i(t) \right].\nn
\end{align}
It follows that, with arbitrary $\Lambda(t)=(\u(t),\r(t),\theta(t),\f(t))$, 
\begin{align}
& \frac{\Ls_i(\bs_i(t))}{\W_i} + \expec [\g_i(\r_i(t);\p_i(t))| \bs_i(t)] \le\frac{\Mone_i(\ks_i)}{\W_i}+ \nn\\
& \frac{\la_i (1-\la_i)\bs_i(t)\ks_i}{\W_i}+  \expec\left[ \frac{\la_i \bs_i(t) \u_i(t)}{\W_i}\!+\!\g_i(\r_i(t);\p_i(t))\vert \bs_i(t) \right].\nn
\end{align}
By summing the above expression over $i=1,\ldots,n$,
\begin{align}
& \sum_{i=1}^n\frac{\Ls_i(\bs_i(t))}{\W_i} + \expec [\g_i(\r_i(t);\p_i(t))| \bs_i(t)] \nn\\
\le & \sum_{i=1}^n\frac{\Mone_i(\ks_i)}{\W_i} \!+\! \frac{\la_i (1-\la_i)\bs_i(t)\ks_i}{\W_i}\!+\!  \expec\big[ \J_{\Pk{3},t}(\Lambda(t)) | \bsn(t)].\nn
\end{align}
where it is clear that minimizing the right hand side of the above inequality over $\Lambda(t)$ is equivalent to minimizing the objective of \Pk{3}. Since $\Lambda^{\mathrm{stat}}(t)$, the disturbance-only stationary policy of \Pk{2}  described in Lemma~\ref{lem:stat_network}, is feasible for \Pk{3}, then the above inequality implies
\begin{align}
& \sum_{i=1}^n\frac{\Ls_i(\bs_i(t))}{\W_i} + \expec [\g_i(\r_i(t);\p_i(t))| \bs_i(t),\Lambda^\mathrm{ol}(t)]\nn\\
\le & \sum_{i=1}^n\frac{\Mone_i(\ks_i)}{\W_i} + \frac{\la_i (1-\la_i)\bs_i(t)\ks_i}{\W_i}+  \expec\big[ \J_{\Pk{3},t}^\star | \bs(t)] \nn \\
\le & \sum_{i=1}^n\!\frac{\Mone_i(\ks_i)}{\W_i}\! +\! \frac{\la_i (1\!-\!\la_i)\bs_i(t)\ks_i}{\W_i}\!+\!  \expec\big[ \J_{\Pk{3},t}(\Lambda^{\mathrm{stat}}(t)) | \bsn(t)]\nn\\
\stackrel{(a)}{=} &\sum_{i=1}^n \frac{\Mone_i(\ks_i)}{\W_i} +\frac{\la_i\bs_i(t)\expec\left[ \ustat_i(t) +(1-\la_i)\ks_i \right]}{\W_i}\nn\\
&+\sum_{i=1}^n\expec [\g_i(\r_i(t);\p_i(t))|\Lambda^{\mathrm{stat}}(t))]\nn \\
\stackrel{(b)}{\le} & \sum_{i=1}^n\frac{\M_i(\ks_i)}{\W_i} +\expec[\g_i(\r_i(t);\p_i(t))|\Lambda^{\mathrm{stat}}(t)]\\
\stackrel{(c)}{\le}& \sum_{i=1}^n\frac{\M_i(\ks_i)}{\W_i} + \Jk{1}^\star.\nn
\end{align}
Here $(a)$ uses the fact that $\ustat(t)$ is induced by a disturbance-only stationary policy;  $(b)$ follows from inequalities
$ |\bs_i(t)|\leq \left(\max\left( (\bmax_i+\ks_i)^2,(\bmin_i+\ks_i)^2\right)\right)^{1/2} $ and
$\left|\expec\left[ \ustat_i(t)\right]+(1-\la_i)\ks_i\right|
\le (1-\la_i) (\max( (\bmax_i+\ks_i)^2,(\bmin_i+\ks_i)^2))^{1/2};
$
and $(c)$ used the following equality $\expec[\sum_{i=1}^n\g_i(\r_i(t);\p_i(t))|\Lambda^{\mathrm{stat}}(t)]  = \Jk{2}^\star$ in Lemma~\ref{lem:stat_network} and $\Jk{2}^\star \le \Jk{1}^\star$.
Taking expectation over $\bsn(t)$ on both sides gives
\begin{align} 
&\expec \left[\sum_{i=1}^n\g_i(\r_i(t);\p_i(t))| \Lambda^{\mathrm{ol}}(t)\right] \label{eq:diff_ly_network}\\
&+ \sum_{i=1}^n\frac{\expec\left[ \L_i(\bs_i(t+1)) -  \L_i( \bs_i(t)) \right]}{\W_i}\nn\leq \sum_{i=1}^n\frac{\M_i(\ks_i)}{\W_i}+ \Jk{1}^\star.\nn
\end{align}
Summing expression \eqref{eq:diff_ly_network} over $t$ from $1$ to $T$, dividing both sides by $ T$, and taking the limit $T\rightarrow \infty$, we obtain the performance bound in expression \eqref{eq:1busiidperf_bdd}.
\end{IEEEproof}


\section{Derivation of the ADMM Algorithm}\label{app:admm}
The first step in deriving the ADMM iterations for the reformulated problem~\eqref{online:ref} is to form the \emph{augmented Lagrangian function} as follows:
\begin{align*}
&L_{\rho}(x, z, \mu, \nu) \\= &\sum_{i=1}^n q_i(x_i) + \mu_i^\intercal (\widehat f_{i, E(i)} - f_{E(i)}) + \frac{\rho}{2}\|\widehat f_{i, E(i)} - f_{E(i)}\|_2^2\\
& + \sum_{e=1}^m h_e(z_e) + \nu_e^\intercal (\widehat \theta_{e, V(e)} - \theta_{V(e)}) + \frac{\rho}{2}\|\widehat \theta_{e, V(e)} - \theta_{V(e)}\|_2^2,
\end{align*}
where $\mu_i \in \real^{|E(i)|}$ and $\nu_e \in \real^{|V(e)|}$ are dual variables for constraints~\eqref{online:ref:c1} and~\eqref{online:ref:c2}, respectively, and $\rho >0$ is a parameter. The centralized ADMM iterates are then
\begin{subequations}
\begin{align}
x^{k+1} &= \argmin_x L_\rho (x, z^k, \mu^k, \nu^k),\label{admm:cen:pn}\\
z^{k+1} &= \argmin_z L_\rho (x^{k+1}, z, \mu^k, \nu^k),\label{admm:cen:pe}\\
\mu^{k+1}_i & = \mu^k_i + \rho(\widehat f_{i, E(i)}^{k+1} - f_{E(i)}^{k+1}), \forall i\in [n],\label{admm:cen:dn}\\
\nu^{k+1}_e &= \nu^k_e + \rho(\widehat \theta_{e, V(e)}^{k+1} - \theta_{V(e)}^{k+1}), \forall e\in [m], \label{admm:cen:de}
\end{align}
\end{subequations}
where $k$ is the iteration count. Let $\eta_i = \mu_i/\rho$ for all $i$ and $\xi_e = \nu_e/\rho$ for all $e$ be the scaled dual variables. Then upon recognizing that updates~\eqref{admm:cen:pn} and~\eqref{admm:cen:dn} are separable across all nodes, and that updates~\eqref{admm:cen:pe} and~\eqref{admm:cen:de} are separable across all edges, we obtain the following distributed ADMM iterates:
\begin{subequations}
\begin{align*}
x^{k+1}_i &= \argmin_{x_i} q_i(x_i) + \frac{\rho}{2} \|\widehat f_{i, E(i)} - f_{E(i)}^k + \eta_i^k\|_2^2\\
&\quad\quad \quad\quad\quad\quad \quad + \sum_{e\in E(i)} \frac{\rho}{2}(\widehat \theta_{e,i}^k- \theta_i + \xi_{e,i}^k)^2,\forall i \in [n],\\ 
z^{k+1}_e &= \argmin_{z_e} h_e(z_e) +\frac{\rho}{2}\|\widehat \theta_{e, V(e)} - \theta_{V(e)}^{k+1} + \xi_e^k\|_2^2\\
&\quad\quad \quad\quad\quad\quad \quad + \sum_{i\in V(e)} \frac{\rho}{2} (\widehat f_{i,e}^{k+1}-f_e +\eta_{i,e}^{k})^2, \forall e \in [m],\\
\eta^{k+1}_i & = \eta^k_i + \widehat f_{i, E(i)}^{k+1} - f_{E(i)}^{k+1}, \forall i\in [n],\\
\xi^{k+1}_e &= \xi^k_e + \widehat \theta_{e, V(e)}^{k+1} - \theta_{V(e)}^{k+1}, \forall e\in [m]. 
\end{align*}
\end{subequations}
The observation that the message passing scheme proposed indeed facilitates the local computation completes this derivation.
\begin{IEEEproof}[Proof of Lemma~\ref{lemma:admmconv}]
Based on the derivation above, it is easy to check that the iterations given above implement the standard two block ADMM with $x$ and $z$ being the (two-block) primal variables, and $(\mu, \nu)$ be the dual variable for the linear equality constraints. The convergence analysis of \cite{BoydADMM} applies directly. The linear convergence rate follows from \eg\, \cite{admmRate} and \cite{2012arXiv1208.3922H}. 

\end{IEEEproof}

\end{document}